\documentclass{article}

\usepackage{arxiv}

\usepackage[utf8]{inputenc} 
\usepackage[T1]{fontenc}    
\usepackage{amsmath}
\usepackage{algorithm}
\usepackage{algorithmic}
\usepackage{hyperref}       
\usepackage{url}            
\usepackage{booktabs}       
\usepackage{amsfonts}       
\usepackage{nicefrac}       
\usepackage{microtype}      
\usepackage{cleveref}       
\usepackage{lipsum}         
\usepackage{graphicx}
\usepackage[numbers]{natbib}
\usepackage{doi}

\title{Adaptive Finite Element Method for Phase Field Fracture Models Based on
Recovery Error Estimates}

\author{
    \href{https://orcid.org/0009-0005-3440-9220}{\includegraphics[scale=0.06]{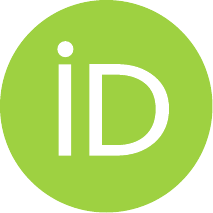}
    \hspace{1mm}Tian Tian}\\
	School of Mathematics and Computational Science\\
	Xiangtan University\\
    China, Xiangtan 411105 \\
    \texttt{tiantian@smail.xtu.edu.cn} \\
	\And
	\href{}{\includegraphics[scale=0.06]{orcid.pdf}
    \hspace{1mm}Chen Chunyu} \\
	School of Mathematics and Computational Science\\
	Xiangtan University\\
    China, Xiangtan 411105 \\
    \texttt{cbtxs@smail.xtu.edu.cn} \\
	\And
	\href{https://orcid.org/0009-0004-7428-9297}{\includegraphics[scale=0.06]{orcid.pdf}
    \hspace{1mm}He Liang} \\
	School of Mathematics and Computational Science\\
	Xiangtan University\\
    China, Xiangtan 411105 \\
    \texttt{lianghesmail@163.com} \\
	\And
	\href{}{\includegraphics[scale=0.06]{orcid.pdf}
    \hspace{1mm}Wei Huayi} 
    \thanks{Corresponding author.}\\
	School of Mathematics and Computational Science, Xiangtan
    University\\
    National Center of Applied Mathematics in Hunan\\
    Hunan Key Laboratory for Computation and Simulation in Science and
    Engineering\\
    China, Xiangtan 411105\\
    \texttt{weihuayi@xtu.edu.cn}\\
}

\renewcommand{\shorttitle}{Adaptive Finite Element Method for Phase Field
Fracture Models Based on Recovery Error Estimates}

\hypersetup{
pdftitle={Adaptive Finite Element Method for Phase Field Fracture Models Based
on Recovery Error Estimates},
pdfsubject={MATHEMATICS},
pdfauthor={Tian Tian, Chen Chunyu, He Liang, Wei Huayi},
pdfkeywords={Adaptive finite method, Fracture, Phase field},
}

\begin{document}
\maketitle

\begin{abstract}
    The phase field model is a widely used mathematical approach for describing
    crack propagation in continuum damage fractures. In the context of phase
    field fracture simulations, adaptive finite element methods (AFEM) are often
    employed to address the mesh size dependency of the model. However, existing
    AFEM approaches for this application frequently rely on heuristic
    adjustments and empirical parameters for mesh refinement. In this paper, we
    introduce an adaptive finite element method based on a recovery type
    posteriori error estimates approach grounded in theoretical analysis. This
    method transforms the gradient of the numerical solution into a smoother
    function space, using the difference between the recovered gradient and the
    original numerical gradient as an error indicator for adaptive mesh
    refinement. This enables the automatic capture of crack propagation
    directions without the need for empirical parameters. We have implemented
    this adaptive method for the Hybrid formulation of the phase field model
    using the open-source software package FEALPy. The accuracy and efficiency
    of the proposed approach are demonstrated through simulations of classical
    2D and 3D brittle fracture examples, validating the robustness and
    effectiveness of our implementation.
\end{abstract}


\section{Introduction}
\label{sec:intro} 
The study of continuous material damage and fracture phenomena holds immense
significance in the realm of materials science and engineering. These phenomena,
influenced by factors such as external loads, thermal stress, and chemical
reactions, can profoundly impact the performance of components and structures,
potentially leading to catastrophic fracture incidents
\cite{anderson2017}.  Understanding and predicting crack propagation
mechanisms is therefore crucial for enhancing the reliability and durability of
materials and structures across diverse engineering sectors, including
aerospace, automotive, construction, and energy
\cite{lemaitre2005, bazant1998, broek2012}. 

The generation of cracks is often accompanied by discontinuities. Addressing
these discontinuities presents significant computational challenges\cite{Newman,
Barenblatt}. The phase field model is used to
describe this problem due to its unique advantages\cite{Kobayashi2010, WU20201}. 
The phase field model describes cracks by introducing continuous phase field
variables, avoiding singularity problems in traditional fracture mechanics and
effectively simulating complex crack paths. It does not require
special treatment of the crack and can directly transform the crack
equation into a simple diffusion equation.
It assumes that the surface tension of the crack naturally exists. It can
easily handle topological changes such as squeezing or merging cracks and is
also easy to combine with various physical equations \cite{Kobayashi2010,
Shen18, MunozReja16, Ruvin23}. The
phase field model can fully consider the microstructure and damage evolution of
materials, thus providing a more accurate description of material behavior. 

Bourdin\cite{Bourdin00} proposed a linear elastic isotropic phase field model in 2000,
wherein both displacement and phase field are treated as linear problem. Amor
\cite{Amor09} introduced a method based on volume-biased strain decomposition in 2009,
which breaks down the total elastic strain energy into volume and deviation
contributions to suppress crack propagation under compression.
Miehe\cite{Miehe10a, Miehe10b}
proposed a spectral decomposition method in 2010 based on crack surface density
energy, facilitating the separation of tension and compression effects on crack
generation and preventing unrealistic cracks under compression. 
To reduce algorithmic complexity, Ambati\cite{Ambati15} introduced a new
approach in 2015. This approach involved a mixed phase field formulation.  The
formulation combined Bourdin's linear elastic isotropic phase field model with
Miehe's maximum historical strain field function.  This combination helped
enhance computational efficiency. Additionally, it prevented crack
reversibility and spurious cracks under compression.

Given the critical role of mesh size in phase field model, especially in
fracture simulation, wherein mesh points at the crack tip must be sufficiently
small, adaptive mesh refinement becomes indispensable.  In 1978,
Babuska\cite{Babuska1978, Babuska1980, Babuska1981} proposed the adaptive finite
element method.  This method can be applied to areas of interest, and numerical
simulations can be performed at a lower cost. In previous studies of fracture
phase field model,  many scholars have proposed different adaptive strategies to
determine mesh refinement regions during fracture evolution\cite{Shanjan2024,
HIRSHIKESH2023109676, XU2023109241}.  Yue et al.\cite{Yue22} utilized local
error estimators based on gradients of the phase-field variable and mechanical
fields to mark elements for refinement, capturing tensile, compressive, and
shear fractures with high precision. Tian et al.\cite{Tian19} combined error
estimation with heuristic marking rules to track changes in the phase-field and
stress fields, making the approach suitable for both quasi-static and dynamic
fracture scenarios. Krishnan et al.\cite{Krishnan22} developed a method that
integrates error estimators with phase-field evolution, marking elements with
rapid phase-field changes for refinement and coarsening those far from the
fracture front to enhance computational efficiency. Assaf et al.\cite{Assaf2022}
adopted an octree-based adaptive method, marking elements for refinement based
on local stress-strain gradients, showing promising results in simulating 3D
brittle fractures.

In this paper, we propose an adaptive finite element method grounded in recovery
type posterior error estimates. This method employs either the average or the
least squares approach to recovery the gradient of the numerical solution of the
phase field from a piecewise constant space to a smoother, piecewise continuous
function space\cite{Huang2012,Zhang2005, LIU20093488, Huang2011}. The disparity between the
smoothed and original numerical solutions serves as an error estimate for
adaptive refinement of mesh elements \cite{Bourdin2014}. 
Compared to the traditional adaptive mesh refinement methods used in phase-field
fracture models, the recovery-type posterior error estimation approach offers
significant advantages. Rather than simply relying on phase-field values,
gradients, or stress-strain variations, it directly evaluates errors through
higher-order approximations based on the concept of error equidistribution,
allowing for precise capture of crack details and automatic identification of
refinement regions, eliminating the need for predefined thresholds.
Additionally, recovery-type posterior error estimation is independent of the
specific problem or finite element method, making it easy to implement with low
computational cost and high robustness, especially in large-scale 3D fracture
simulations.

We leverage the open-source numerical PDE solver FEALPy\cite{Wei2017}, utilizing
object-oriented and array-based programming techniques to efficiently implement
the proposed algorithms. FEALPy is developed entirely using standard Python
scientific computing libraries, such as NumPy \cite{NumPy}, SciPy \cite{SciPy},
and Matplotlib \cite{Matplotlib}, and
provides a rich set of mesh data structures and adaptive algorithms, including
commonly used 1D, 2D, and 3D mesh objects, as well as refinement methods like
bisection and red-green refinement. Our implementation employs array-based
programming to maintain code simplicity while ensuring high computational
efficiency. Additionally, we introduce matrix assembly techniques that avoid
numerical integration, and GPU acceleration to solve the discrete systems,
significantly enhancing algorithm performance. Numerous classical examples are
presented to validate the efficiency of our implementation.

The remainder of this paper is organized as follows. In section
\ref{sec:model}, we introduce the energy function, control equations of the
fracture phase field model, and the formulation phase field model. In
section \ref{sec:Algorithm}, we introduce the algorithm used in this article,
including the adaptive finite element method, posterior error estimates, and the
corresponding program algorithm flow. In section \ref{sec:application}, we
present the numerical experimental results of several classic examples. In
section \ref{sec:summary}, we have summarized the the contributions and findings
of this study.

\section{Mathematical model}
\label{sec:model}
\subsection{Crack energy function}
Given a material region $\Omega$, if we apply a force $\boldsymbol{f}$, it will
lead to deformation of the material. This deformation can be measured by
displacement $\boldsymbol{u}$. If cracks occur during this deformation process,
the material can be classified into two states, and we can use phase function $d$ to
describe them: cracked ($d=1$) and uncracked ($d=0$)\cite{Bourdin2014,
LI201944}. At this point, fracture energy density function can be expressed as

\begin{equation}
\gamma(d, \nabla d) = \frac{1}{c_h}\left(\frac{h(d)}{l_0} + l_0|\nabla d|^2\right)
\end{equation}

\begin{itemize}
  \item $l_0 > 0$: a small regularization length scale parameter
      \cite{Bourdin2014, Irwin48, Irwin57}, used to control the width of the smoothly approximated crack.
  \item $h(d)$: a continuously and strictly increasing function.
  \item Normalization constant $c_h = \int_0^1\sqrt{h(d)}$: a constant related to $h(d)$ that ensures the regularization crack surface converges to the shaped crack surface as $l_0 \to 0$.
\end{itemize}

There are generally two options for fracture energy density function, namely:
\begin{itemize}
  \item AT1 model: $h(d)=d, c_h=\frac{8}{3}$ (finite width damage band, including pure elastic stage)
  \item AT2 model: $h(d)=d^2, c_h=2$ (infinite width damage band)
\end{itemize}
In this article, we use the AT2 model. The formulation presented here is based
on the format proposed by Miehe\cite{Miehe10a, Miehe10b} et al. The expressions
for fracture energy $E_c(d)$ are defined as follows:
\begin{equation}
\begin{aligned}
E_c(d) &= G_c\int_{\Omega}\gamma(d, \nabla d)\mathrm{d} \boldsymbol{x}
\end{aligned}
\end{equation}
Here, $G_c > 0$ denotes the critical energy release rate of the material\cite{Griffith}.

\subsection{Strain energy function with crack propagation}
The energy stored by materials under external forces is called strain energy.
For elastic materials, when external forces are eliminated, this energy will be
completely released. However, in the case of crack propagation during the
deformation process, a portion of the stored strain energy will be lost.
Therefore, it is necessary to introduce the energy degradation function $g(d)$
in the expression of strain energy to consider crack propagation. The energy
degradation function is defined as 
\begin{equation}
g(d) = (1-d)^2 + \epsilon,
\end{equation}
where $0 \le \epsilon \ll 1$ is introduced to ensure numerical stability. Of
course, the energy degradation function can also choose other
formats\cite{Borden2012, Karma2001, Sargado2017, wu2018}. At this point, strain energy
is expressed as
\begin{equation}
    E_s(\boldsymbol{u}, d) = \int_{\Omega}
    g(d)e_s(\boldsymbol{\varepsilon})\mathrm{d}\boldsymbol{x}-W(\boldsymbol{u}),
\end{equation}
where, $W(\boldsymbol{u})$ is the external force. The strain tensor
$\boldsymbol{\varepsilon}$ is given by:
\begin{equation}
\boldsymbol{\varepsilon} = \frac{1}{2}(\boldsymbol{u} \otimes \nabla + \nabla
\otimes \boldsymbol{u}).
\end{equation}

However, the above expression does not distinguish between tensile and compressive stress
states, and spurious cracks may form in regions under compression. To avoid this
situation, the commonly used processing method invloves dividing the energy density
$e_s(\boldsymbol{\varepsilon})$ is divided into positive (related to tension
state) and negative (compression state) parts, and applying the degradation
function only to the positive component. This results in the following
expression: 

\begin{equation}
    E_s(\boldsymbol{u}, d) =\int_{\Omega} \left[g(d) e_s^+(\boldsymbol{\varepsilon}) +
    e_s^-(\boldsymbol{\varepsilon})\right]\mathrm{d}\boldsymbol{x}-W(\boldsymbol{u}).
\end{equation}

\subsection{Governing equations}

According to the previous energy formula, we can obtain the total potential energy function:
\begin{equation}
\begin{aligned}
L &= K(\dot{\boldsymbol{u}}) - \Pi(\boldsymbol{u}, d)\\
&=
\frac{1}{2}\int_{\Omega}\rho\dot{\boldsymbol{u}}\cdot\dot{\boldsymbol{u}}\,\mathrm{d}\boldsymbol{x}
- \Pi(\boldsymbol{u}, d),
\end{aligned}
\end{equation}
where
\begin{equation}
\begin{aligned}
\Pi(\boldsymbol{u}, d) &= E_s(\boldsymbol{u},d) + E_c(d)\\
&= \int_{\Omega} [(1-d)^2 + \epsilon] e_s^+(\boldsymbol{\varepsilon}) + e_s^-(\boldsymbol{\varepsilon})\mathrm{d}\boldsymbol{x} \\
&\quad + \int_{\Omega}\frac{G_c}{2}\left[\frac{1}{l_0}d^2 + l_0|\nabla
d|^2\right]\mathrm{d}\boldsymbol{x} - \int_{\Omega}\boldsymbol{f}\cdot
\boldsymbol{u}\mathrm{d}\boldsymbol{x}.
\end{aligned}
\end{equation}

According to the Hamiltonian principle, the variation of the Lagrange operator
can lead to the Euler Lagrange equation:
\begin{equation}
\begin{aligned}
\begin{cases}
\rho \ddot{\boldsymbol{u}} - \text{div}\boldsymbol{\sigma} = \boldsymbol{f}\\
g'(d) e_s^+(\boldsymbol{\varepsilon}) + G_c\gamma'(d, \nabla d) = 0.
\end{cases}
\end{aligned}
\end{equation}
The stress tensor $\boldsymbol \sigma$ is:
\begin{equation}
\boldsymbol{\sigma}(\boldsymbol{\varepsilon}) := g(d)
\frac{\partial e_s^+(\boldsymbol{\varepsilon})}{\partial
\boldsymbol{\varepsilon}} + \frac{\partial
e_s^-(\boldsymbol{\varepsilon})}{\partial \boldsymbol{\varepsilon}}.
\end{equation}

\subsection{Hybrid phase field model}

For the "Isotropic" Model \cite{Bourdin00}, strain energy density decomposition is defined as:
\begin{equation}
e_s^+(\boldsymbol{\varepsilon}) = e_s(\boldsymbol{\varepsilon}), \qquad
e_s^-(\boldsymbol{\varepsilon}) = 0.
\end{equation}
Here, the strain energy density of linear isotropic materials is defined as:
\begin{equation}
e_s = \frac{\lambda}{2} \text{tr}(\boldsymbol{\varepsilon})^2 + \mu
\boldsymbol{\varepsilon} : \boldsymbol{\varepsilon}.
\end{equation}
Where, $\lambda$ is the  Lamé's first parameter, $\mu$ is the Lamé's second parameter
(shear modulus), and $\mu > 0$, $\lambda + 2\mu > 0$. Their relationships with
Young's modulus $E$, Poisson's ratio $\nu$ are as follows:
\begin{equation}
    \lambda = \frac{\nu E}{(1 + \nu)(1 - 2\nu)} , \qquad
    \mu = \frac{E}{2(1 + \nu)} .
\end{equation}

Miehe \cite{Miehe10a, Miehe10b} proposed a spectral decomposition method, where strain energy density decomposition is defined as:
\begin{equation}
\begin{cases}
e_s^+(\boldsymbol{\varepsilon}) &= \frac{\lambda}{2}<\text{tr}(\boldsymbol{\varepsilon})>_+^2 + \mu \text{tr}(\boldsymbol{\varepsilon}_+^2) \\
e_s^-(\boldsymbol{\varepsilon}) &=
\frac{\lambda}{2}<\text{tr}(\boldsymbol{\varepsilon})>_-^2 + \mu
\text{tr}(\boldsymbol{\varepsilon}_-^2),
\end{cases}
\end{equation}
and,
\begin{equation}
\boldsymbol{\varepsilon}_{\pm} =
\sum_{i=0}^{n-1}<\varepsilon_i>_{\pm}\boldsymbol{n}_i\otimes\boldsymbol{n}_i,
\end{equation}
where,
\begin{itemize}
  \item $\varepsilon_i$ represents eigenvalues.
  \item $\boldsymbol{n}_i\otimes \boldsymbol{n}_i$ represents eigenvectors.
  \item $<\cdot>_{\pm}$ is Macqualay operation, which is defined as: $<x>_{\pm} = \frac{1}{2}(x \pm |x|).$
\end{itemize}

As crack growth is an irreversible phenomenon, Miehe proposed a maximum strain field function:
\begin{equation}
\mathcal{H}(\boldsymbol{x}, t) =
\max_{s\in[0,t]}e_s^+(\boldsymbol{\varepsilon}(\boldsymbol{x},s))
\end{equation}
to replace $e_s^+(\boldsymbol{\varepsilon})$ in the equation
$g'(d)e_s^+(\boldsymbol{\varepsilon}) + G_c\gamma'(d, \nabla d) = 0$, which can
prevent the occurrence of non-physical crack closure.

In the hybrid mixed format \cite{Ambati15}, the "Isotropic" model is combined
with Miehe's\cite{Miehe10a} maximum historical strain field function to obtain a linear mixed model, where the expression for the stress tensor is
\begin{equation}
\boldsymbol{\sigma}(\boldsymbol{\varepsilon}) = g(d)\frac{\partial
e_s(\boldsymbol{\varepsilon})}{\partial \boldsymbol{\varepsilon}} = g(d)[\lambda
(\text{tr}(\boldsymbol{\varepsilon})\boldsymbol{I}) + 2\mu
\boldsymbol{\varepsilon}].
\end{equation}
And the fourth-order elastic tangent operator tensor:
\begin{equation}
\mathbb{C}:= \frac{\partial \boldsymbol{\sigma}}{\partial
\boldsymbol{\varepsilon}} = g(d)(\lambda \boldsymbol{I} \otimes \boldsymbol{I} + 2\mu
\mathbb{I}).
\end{equation}

Therefore, the final governing equations become:
\begin{equation}
\begin{aligned}
\begin{cases}
\rho \ddot{\boldsymbol{u}} - \text{div}\boldsymbol{\sigma} = \boldsymbol f \\
g'(d) \mathcal{H} + G_c\gamma'(d, \nabla d) = 0,
\end{cases}
\end{aligned}
\end{equation}
with
\begin{equation}
    \begin{cases}
        d=0,
        \text{ on }\partial \Omega\\
        \boldsymbol\sigma \cdot \boldsymbol n=\boldsymbol g,
        \text{ on } \Gamma_N.
    \end{cases}
\end{equation}

\section{Algorithm design}
\label{sec:Algorithm}
\subsection{Nonlinear finite element discretization}
In this article, we only consider quasi-static fracture, so we can abandon term
$\rho \ddot{\boldsymbol u}$. Using the finite element method, the variational
formulation of the governing equations can be derived, yielding the weak form as
follows:
\begin{equation}
\begin{aligned}
\begin{cases}
    (\boldsymbol\sigma(\boldsymbol{u}), \boldsymbol{\varepsilon}(\boldsymbol{v}) )_{\Omega} =
    (\boldsymbol{f},\boldsymbol{v} )_{\Omega}  + <\boldsymbol g,
    \boldsymbol{v}>_{\Gamma_N}, \qquad 
    &\boldsymbol{v} \in (H^1)^{dim}(\Omega)\\
    -2((1-d)\mathcal H, \omega)_{\Omega} + G_c l_0(\nabla d, \nabla
    \omega)_{\Omega} + \frac{G_c}{l_0}(d, \omega)_{\Omega} =0, \qquad
    &\omega\in H^1(\Omega),
\end{cases}
\end{aligned}
\end{equation}
where, $dim$ is the dimension of the problem.

In this article, we use triangle and tetrahedral mesh for discretization. We
discretize $H ^ 1(\Omega)$ by finite element space, and the maximum strain field
$\mathcal H \in L^2(\Omega)$. The discretization weak form is:
\begin{equation}
    \begin{aligned}
        \begin{cases}
            (\boldsymbol\sigma(\boldsymbol{u}_h),
            \boldsymbol{\varepsilon}(\boldsymbol{v_h}) )_{\Omega} =
            (\boldsymbol{f},\boldsymbol{v_h} )_{\Omega} + <\boldsymbol g,
            \boldsymbol v_h>_{\Gamma_N} \\
            -2((1-d_h)\mathcal H, \omega_h)_{\Omega} + G_c l_0(\nabla d_h, \nabla
            \omega_h)_{\Omega} + \frac{G_c}{l_0}(d_h, \omega_h)_{\Omega} =0.
        \end{cases}
    \end{aligned}
\end{equation}

Let the basis functions for the Lagrange scalar space be denoted as $\boldsymbol
\phi$, and the basis functions for the Lagrange vector space be denoted as
$\boldsymbol \Phi$. 
At this point, the discretization weak form can be obtained:
\begin{equation}
\begin{aligned}
\begin{cases}
    (\boldsymbol\sigma(\boldsymbol{u}_h),
    \boldsymbol{\varepsilon}(\boldsymbol{\Phi}) )_{\Omega} =
    (\boldsymbol{f},\boldsymbol{\Phi} )_{\Omega} + <\boldsymbol g, \boldsymbol
    \Phi>_{\Gamma_N}\\
-2((1-d_h)\mathcal H, \phi)_{\Omega} + G_c l_0(\nabla d_h, \nabla
\phi)_{\Omega} + \frac{G_c}{l_0}(d_h, \phi)_{\Omega} =0.
\end{cases}
\end{aligned}
\end{equation}

The residual vectors $\boldsymbol R_0$ and $\boldsymbol R_1$ are defined as:
\begin{equation}
\begin{aligned}
        \boldsymbol{R}_0 &= (\boldsymbol f, \boldsymbol{\Phi} )_{\Omega}
        + <\boldsymbol g, \boldsymbol \Phi>_{\Gamma_N}
        - (\boldsymbol{\sigma}(\boldsymbol u_h), \boldsymbol
        \varepsilon(\boldsymbol{\Phi}))_{\Omega} \\
        \boldsymbol{R}_1 &= 2((1-d_h)\mathcal H(\boldsymbol u_h), \boldsymbol\phi)_{\Omega} - G_c
        l_0(\nabla d_h, \nabla \boldsymbol\phi)_{\Omega} - \frac{G_c}{l_0}(d_h,
        \boldsymbol\phi)_{\Omega}.
\end{aligned}
\end{equation}

We employ the Newton-Raphson method for iterative solving, which leads to the
following system of equations:
\begin{equation}
\begin{aligned}
\begin{bmatrix}
    -\frac{\partial \boldsymbol R_0}{\partial \boldsymbol{u_h}} &
    -\frac{\partial \boldsymbol R_0}{\partial d_h} \\
    -\frac{\partial \boldsymbol R_1}{\partial \boldsymbol{u_h}} &
    -\frac{\partial \boldsymbol R_1}{\partial d_h}
\end{bmatrix}
\begin{bmatrix}
\Delta \boldsymbol{u_h} \\ 
\Delta d_h
\end{bmatrix}
=
\begin{bmatrix}
\boldsymbol R_0 \\ 
\boldsymbol R_1
\end{bmatrix}
\end{aligned}
\end{equation}
where, $\Delta \boldsymbol{u_h} = \boldsymbol{u_h}^{k} - \boldsymbol{u_h}^{k-1}$,
$\Delta d_h = d_h^{k} - d_h^{k-1}$, and $\boldsymbol{u_h}^{k}$, $d_h^{k}$ are
the solutions at the $k$-th iteration. Thus, the element stiffness matrix can
be:

\begin{equation}
\begin{aligned}
    -\frac{\partial {\boldsymbol{R}}_0}
    {\partial \boldsymbol{u}_h} &= \frac{\partial}{\partial
    \boldsymbol{u}_h}(\boldsymbol\sigma(\boldsymbol u_h), \boldsymbol{\varepsilon}(\boldsymbol{\Phi})
    )_{\Omega} \\
    -\frac{\partial {\boldsymbol{R}}_0}{\partial d_h} &= \frac{\partial
        \boldsymbol{\sigma}}{\partial d_h}(\nabla\boldsymbol{\Phi}^{T}, \boldsymbol\phi)_{\Omega} \\
    -\frac{\partial \boldsymbol{R}_1}
    {\partial \boldsymbol{u}_h} &= -2(1-d_h)\frac{\partial \mathcal{H}}{\partial
        \boldsymbol{\varepsilon}} (\boldsymbol\phi^{T}, \nabla \boldsymbol{\Phi})_{\Omega}\\
    -\frac{\partial \boldsymbol{R}_1}
    {\partial d_h} &= 2(\boldsymbol\phi^T\mathcal H, \boldsymbol\phi)_{\Omega} + G_c
    l_0(\nabla \boldsymbol\phi^T, \nabla \boldsymbol\phi)_{\Omega} +
    \frac{G_c}{l_0}(\boldsymbol\phi^T, \boldsymbol\phi)_{\Omega},
\end{aligned}
\end{equation}
where the derivative of the stress tensor with respect to the displacement:
\[
    \frac{\partial}{\partial \boldsymbol{u_h}}(\boldsymbol\sigma(\boldsymbol u_h),
    \boldsymbol{\varepsilon}(\boldsymbol{\Phi}) )_{\Omega} = \frac{\partial
    \boldsymbol{\sigma}}{\partial \boldsymbol{\varepsilon}}
    (\boldsymbol{\varepsilon}(\boldsymbol{\Phi}),
    \boldsymbol{\varepsilon}(\boldsymbol{\Phi}))_{\Omega}.
\]

Because the total energy is a non convex function of the overall unknown
quantity $(\boldsymbol u, d)$, but it is a separate convex function of
displacement $\boldsymbol u$ or phase field $d$.
Therefore, in a robust Staggered strategy, $\frac{\partial \boldsymbol{R}_1}{\partial
\boldsymbol{u_h}}, \frac{\partial \boldsymbol{R}_0}{\partial d_h}$ can be
ignored, to independently solve for the displacement variable $\boldsymbol{u_h}$
and the phase field variable $d_h$. This independence simplifies the problem and
allows for separate updates of these variables without direct coupling in every
iteration or time step. In this article, we use the linear element method, so we
can use the method of no numerical integration for matrix assembly, including
$(\boldsymbol \phi, \boldsymbol \phi)_{\Omega}$. Here, we will not elaborate
further. The specific content can be found in the FEALPy software
package\cite{Wei2017}.

\subsection{An adaptive method based on recovery type posterior error estimates}

In this paper, we use a recovery-type posterior error estimates method for
error calculation, which is based on post-processing techniques. The commonly
used methods are: weighted average, least squares, local or global projection,
etc. The post-processed solution is used instead of the gradient of the true
solution and the numerical solution for error estimates\cite{Huang2012,
Huang2011}. Based on the
characteristics of the phase field function, it can be inferred that the
distribution of the phase field can effectively capture the occurrence of
cracks. Therefore, we use the phase field function for posterior error
estimates.

The process of the adaptive mesh refinement algorithm based on posterior error
estimates is as follows: after the numerical solution is calculated on the
current mesh, the element error and global error are estimated by using the
posterior error estimates. The cells that need to be refined are labeled
according to the calculated element error and the selected labeling strategy, and
finally, the selected refinement method is used to refine the mesh. 

Below is an introduction to commonly used gradient recovery methods. The simple
averaging method is:
\[
    \mathcal {R}_h d_h(\boldsymbol{x}) = \frac{1}{m_{\boldsymbol{x}}} \sum_{\tau \in
\Gamma_{\boldsymbol{x}}} \nabla d_h(\boldsymbol{x})|_{\tau}.
\]
Here, $m_{\boldsymbol{x}}$ is the number of elements with node $\boldsymbol{x}$ as
the vertex, and $\Gamma_{\boldsymbol{x}}$ is the set of elements with node
$\boldsymbol{x}$ as the vertex.

Area averaging method is:
\[
\mathcal R_h d_h(\boldsymbol{x}) = \frac{\sum_{\tau \in \Gamma_{\boldsymbol{x}}} |\tau|
\nabla d_h(\boldsymbol{x})|_{\tau}}{\sum_{\tau \in \Gamma_{\boldsymbol{x}}}
|\tau|}.
\]
Here, $|\tau|$ represents the area of element $\tau$.

Harmonic area averaging method is:
\[
\mathcal R_h d_h(\boldsymbol{x}) = \frac{\sum_{\tau \in \Gamma_{\boldsymbol{x}}}
\frac{1}{|\tau|} \nabla d_h(\boldsymbol{x})|_{\tau}}{\sum_{\tau \in
\Gamma_{\boldsymbol{x}}} \frac{1}{|\tau|}}.
\]

Angle averaging method is:
\[
\mathcal R_h d_h(\boldsymbol{x}) = \frac{\sum_{\tau \in \Gamma_{\boldsymbol{x}}}
\alpha_\tau \nabla d_h(\boldsymbol{x})|_{\tau}}{\sum_{\tau \in
\Gamma_{\boldsymbol{x}}} \alpha_\tau}.
\]
Here, $\alpha_\tau$ represents the degree of the angle of element $\tau$ with a vertex of $\boldsymbol{x}$.

Distance average is:
\[
\mathcal R_h d_h(\boldsymbol{x}) = \frac{\sum_{\tau \in \Gamma_{\boldsymbol{x}}} l_\tau
\nabla d_h(\boldsymbol{x})|_{\tau}}{\sum_{\tau \in \Gamma_{\boldsymbol{x}}}
l_\tau}.
\]
Here, $l_\tau$ represents the distance from the center of gravity of element $\tau$ to $\boldsymbol{x}$.

On the element $\tau$, let $\mathcal R_h$ be the gradient reconstruction operator, and for the phase field $d$ there is:
\[
\eta_{\tau} = \|\mathcal R_h d_h - \nabla d_h \|_{0,\tau}.
\]

The above equation can be used to estimate the error $\| \nabla e_h \|_{\Omega,
\tau}$. If $\mathcal R_h$ is super-convergent \cite{Zhang2005}, then
\[
\|\mathcal R_h d_h - \nabla d_h \|_{0,\Omega} = o(\| \nabla d - \nabla d_h \|_{0, \Omega}),
\]
therefore,
\[
\frac{\|\mathcal R_h d_h - \nabla d_h \|_{0,\Omega}}{\| \nabla d - \nabla d_h
\|_{0, \Omega}} = \frac{\| (\mathcal R_h d_h - \nabla d) + ( \nabla d - \nabla
d_h) \|_{0,\Omega}}{\| \nabla d - \nabla d_h \|_{0, \Omega}} \to 1.
\]

The posterior error estimate obtained in this case is asymptotically
accurate. Based on the calculated posterior error estimates, the error on each
element can be obtained. Units with larger errors are selected for refinement, and the error
size is judged according to the labeling strategy to determine whether it meets
the standard. Detailed theories and application examples can be found in papers
such as \cite{zienkiewicz1992part1, zienkiewicz1992part2, ainsworth2000, dai2004}. 
There are two commonly used labeling strategies, where $\eta$ is
the total estimates error, $\eta_\tau$ is the estimates error on element $\tau$,
and $0 < \theta < 1$ is the labeling parameter.

\begin{itemize}
    \item The Maximum Criterion method: Make $\eta_{\max} = \max_{\tau \in
        \Gamma} \eta_{\tau}$, mark the elements of $\eta_\tau > \theta \eta_{\max}$.
    \item The $L^2$ Criterion method: Sort the element estimates errors and
        place the corresponding elements in the labeled element set $\mathcal{M}$ from
        largest to smallest, until $\sum_{\tau \in \mathcal{M}} \eta_\tau^2 >
        \theta \eta^2$ is satisfied.
\end{itemize}

\subsection{Overall Algorithm Process}
\begin{algorithm}[!ht]
\caption{Solve Process}
	\begin{algorithmic}[1]
        \FOR {the iterations $k$ and the current time step $n$}
        \STATE calculate the displacement stiffness matrix
        $-\left(\frac{\partial \boldsymbol{R}_0}{\partial
        \boldsymbol{u_h}}\right)^k$ and the right-hand side term
        $\boldsymbol{R}_0^k$.  
        \STATE solve for displacement
        $\boldsymbol{u_h}^{k+1}$: \[ -\left(\frac{\partial
            \boldsymbol{R}_0}{\partial \boldsymbol{u_h}}\right)^k
            \Delta \boldsymbol{u_h}^k = \boldsymbol{R}_0^k, \quad
            \boldsymbol{u_h}^{k+1} = \boldsymbol{u_h}^k + \Delta
        \boldsymbol{u_h}^k, \]
        \STATE substitute $\boldsymbol{u_h}^{k+1}$ into the $\mathcal{H}$ function
        to obtain: \[ \mathcal{H}^{k+1}(\boldsymbol{u_h}^{k+1}) = \max
        e^+(\boldsymbol{\varepsilon}(\boldsymbol{u_h}^{k+1})), \] 
        \STATE using the
        updated $\mathcal{H}$ function, calculate the phase field stiffness
        matrix $-\left(\frac{\partial \boldsymbol{R}_1}{\partial d_h}\right)^k$
        and the right-hand side term $\boldsymbol{R}_1^k$.  
        \STATE solve for
        the phase field value $d_h^{k+1}$: \[ -\left(\frac{\partial
            \boldsymbol{R}_1}{\partial d_h}\right)^k \Delta d_h^k =
    \boldsymbol{R}_1^k, \quad d_h^{k+1} = d_h^k + \Delta d_h^k.  \] 
        \STATE calculate the posterior error estimate $\eta$, adaptive mesh refinement.  
        \STATE compute the relative residual
    $error=
    \max\left\{\frac{\|\boldsymbol{R}_0^{(n+1)}\|}{\|\boldsymbol{R}_0^{(0)}\|},
    \frac{\|\boldsymbol{R}_1^{(n+1)}\|}{\|\boldsymbol{R}_1^{(0)}\|}\right\}$.
        \IF {error < $1e-5$ }
            \STATE $k+1$.
        \ELSE 
            \STATE $n+1$.  
        \ENDIF
    \ENDFOR
    \end{algorithmic}
\end{algorithm}

\begin{figure}[htbp]
    \centering
    \includegraphics[width=\textwidth]{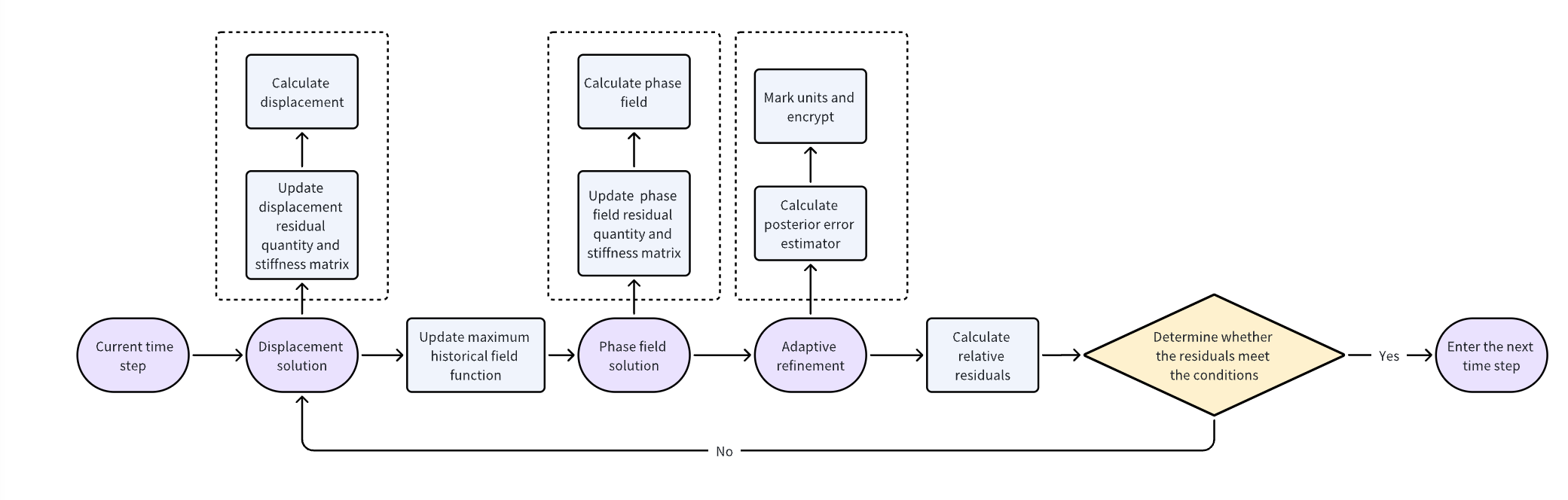}
    \caption{Algorithm Flow}
    \label{fig:alg}
\end{figure}

Based on FEALPy\cite{Wei2017}, we have developed the program in a modular fashion, employing a
unified interface to significantly enhance the flexibility and efficiency of
experimental testing. Figure \ref{fig:alg} is the overall algorithm
implementation flowchart. Each algorithm's program module is independent,
allowing for  the selection of appropriate calculation methods. In our program
design, the choice of solution method is dimension-agnostic, allowing
independent selection of adaptive mesh refinement for both 2D and 3D
calculations.

To optimize the algorithm, we have introduced a matrix calculation method
without numerical integration assembly in the matrix assembly module.
Additionally, in solving algebraic equation systems, we have integrated various
methods including direct solvers, iterative solvers \cite{SciPy}, and GPU iterative
solvers \cite{CuPy}. It allows for the selection of appropriate methods, resulting in
significant computational efficiency gains. For Mesh refinement methods, bisect
refinement and triangle red-green refinement based on half-edge-mesh can be
selected according to the characteristics of the mesh.

The simulations are performed on a machine with the following specifications:

\begin{itemize}
    \item \textbf{Memory:} 64.0 GiB
    \item \textbf{Processor:} Intel\textsuperscript{\textregistered} Xeon(R) Gold 5118 CPU @ 2.30GHz $\times$ 48 cores
    \item \textbf{Graphics:} NVIDIA Corporation GP107GL [Quadro P620]
\end{itemize}

\begin{table}[!ht]
    \centering
    \caption{Comparison of Solve Times with and without GPU for a 3D Model}
    \begin{tabular}{|c|c|c|c|}
        \hline
        & GPU (Iterative)  & CPU (Iterative) & CPU (Direct) \\
        \hline
        Displacement Solve Time (seconds)& 4.90454 & 14.18117 & 12391.396 \\
        \hline
        Phase Field Solve Time (seconds) & 0.01722 & 0.07701 & 64.9935 \\
        \hline
    \end{tabular}
    \label{table:gpu_comparison}
\end{table}

The table (\ref{table:gpu_comparison}) compares the computational performance
for solving a 3D model with 56,729 nodes, using GPU acceleration versus CPU
processing. It presents the solve times for both displacement and phase field
calculations under different methods.

Displacement solve time: Using GPU, the time required to solve displacement is
4.90454 seconds. In comparison, using the CPU with the iterative method takes
14.18117 seconds, and the direct method takes 12391.396 seconds. The GPU
solution is approximately 65.5\% faster than the iterative method and about
99.96\% faster than the direct method.

Phase field solve time: For phase field calculations, the GPU solution takes
0.01722 seconds. The CPU iterative method requires 0.07701 seconds, while the
direct method takes 64.9935 seconds. The GPU solution is approximately 77.6\%
faster than the iterative method and about 99.97\% faster than the direct method.

Overall, the use of GPU acceleration 
results in substantial time savings and improved computational efficiency,
making these techniques highly advantageous in practical applications,
especially for solving large-scale 3D models.

\section{Algorithm validation and application}
\label{sec:application}
We test several classical examples, namely: a model with a rigid
circular inclusion embedded in a square plate, a model with a unilateral notch
in which the plate is subjected to tensile and shear forces, an L-shaped plate
with mixed tensile and compressive loads, and a 3D model with a single slice
under tensile force.

\subsection{Square Model with a Circular Notch: Upward Stretch on the Upper Boundary}

\begin{figure}[htbp]
    \centering
    \includegraphics[width=0.4\textwidth]{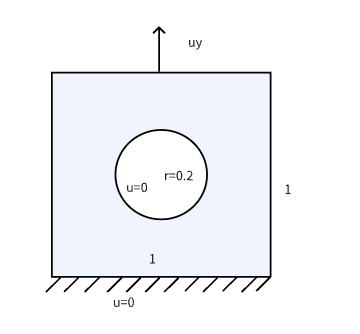}
    \caption{Geometry and boundary conditions of the square model with a
    circular notch (the unit is mm)}
    \label{fig:geometry0}
\end{figure}

A rigid circular inclusion is embedded in a square plate, and a vertical
displacement is applied to its top surface.  The computational domain $\Omega$
is a rectangular region $[0, 1] \times [0, 1]$ with a circular cutout centered
at $(0.5, 0.5)$, having a radius of $0.2$, Figure \ref{fig:geometry0} shows the
geometric shape of the model.For the phase field on the boundaries, zero
Dirichlet boundary conditions are applied. Parameters are given as follows:
$G_c=1$, $l_0=0.02$, $E=200$, $\nu=0.2$\cite{Bourdin00}.

Boundary conditions are as follows:
\begin{itemize}
    \item Dirichlet boundary at the upper boundary $y=1$, with displacement
        increments as follows for the first 5 steps: $\Delta
        u_y=1.4\times10^{-2} \text{ mm}$, and for the subsequent 25 steps:
        $\Delta u_y=2.2\times10^{-3} \text{ mm}$;
    \item Internal boundary: Displacement is 0 along the boundary of the
        circular hole.
\end{itemize}

\begin{figure}[htbp]
    \centering
    \includegraphics[width=\textwidth]{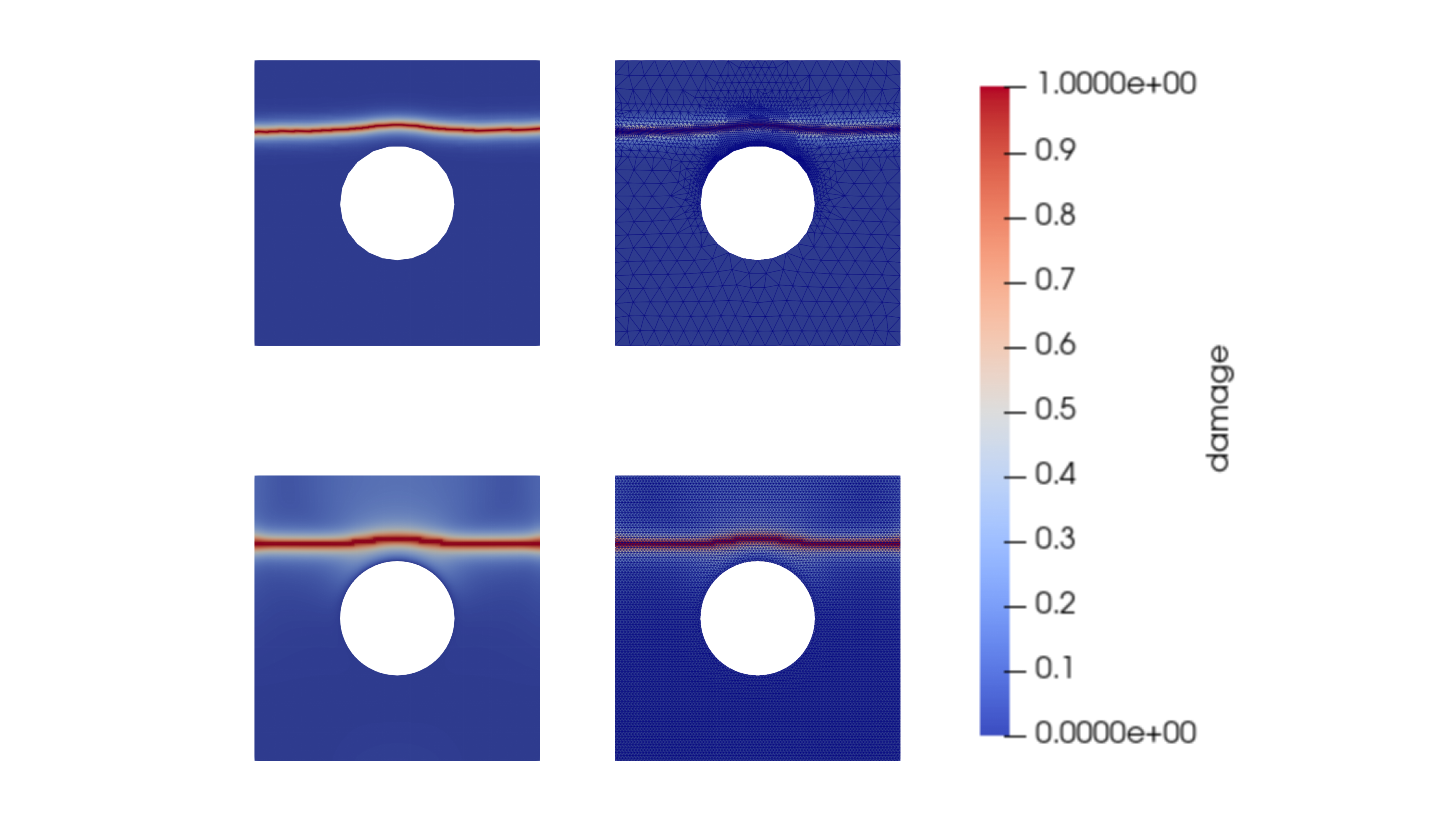}
    \caption{The final results of the square model with a circular notch. The
        upper figures show the results of adaptive mesh refinement (initial
        mesh size h=0.05), and the lower figures show the results without
        using adaptive (mesh size h=0.01). The left figures show the final phase
        field values, and the right figures show the final mesh situation.}
    \label{fig:mesh0_result}
\end{figure}

As shown in Figure \ref{fig:mesh0_result}, the adaptive finite element
simulation displays the final phase field values, crack states, and the
resulting mesh refinement. With an initial mesh size of 0.05, the mesh is
automatically refined in regions that require higher accuracy during the
simulation. This dynamic refinement results in an increase in the number of mesh
elements from 640 to 9558, ensuring high simulation accuracy. The adaptive
refinement approach significantly improves computational efficiency, completing
the entire simulation in only 86 seconds.

In contrast, if the adaptive refinement is not used, we must ensure that the
mesh size $h<l_0/2$. At this poin, the non-adaptive simulation utilizes a fine
mesh size of 0.01 from the start and maintains a constant mesh element of 19224
throughout the simulation. While this guarantees uniform high accuracy
throughout the simulation domain, it also results in wasted computational
resources in regions where finer meshes are not necessary. Consequently, the
simulation time is significantly longer, taking 695 seconds to complete.

Despite the differences in mesh refinement and simulation time, the results from
both the adaptive and non-adaptive simulations exhibit a high degree of
consistency in terms of the computed phase field values and crack states.
However, the adaptive refinement method stands out as a preferred approach due
to its ability to balance computational accuracy and efficiency through dynamic
mesh refinement.

\begin{figure}[htbp]
    \centering
    \includegraphics[width=0.5\textwidth]{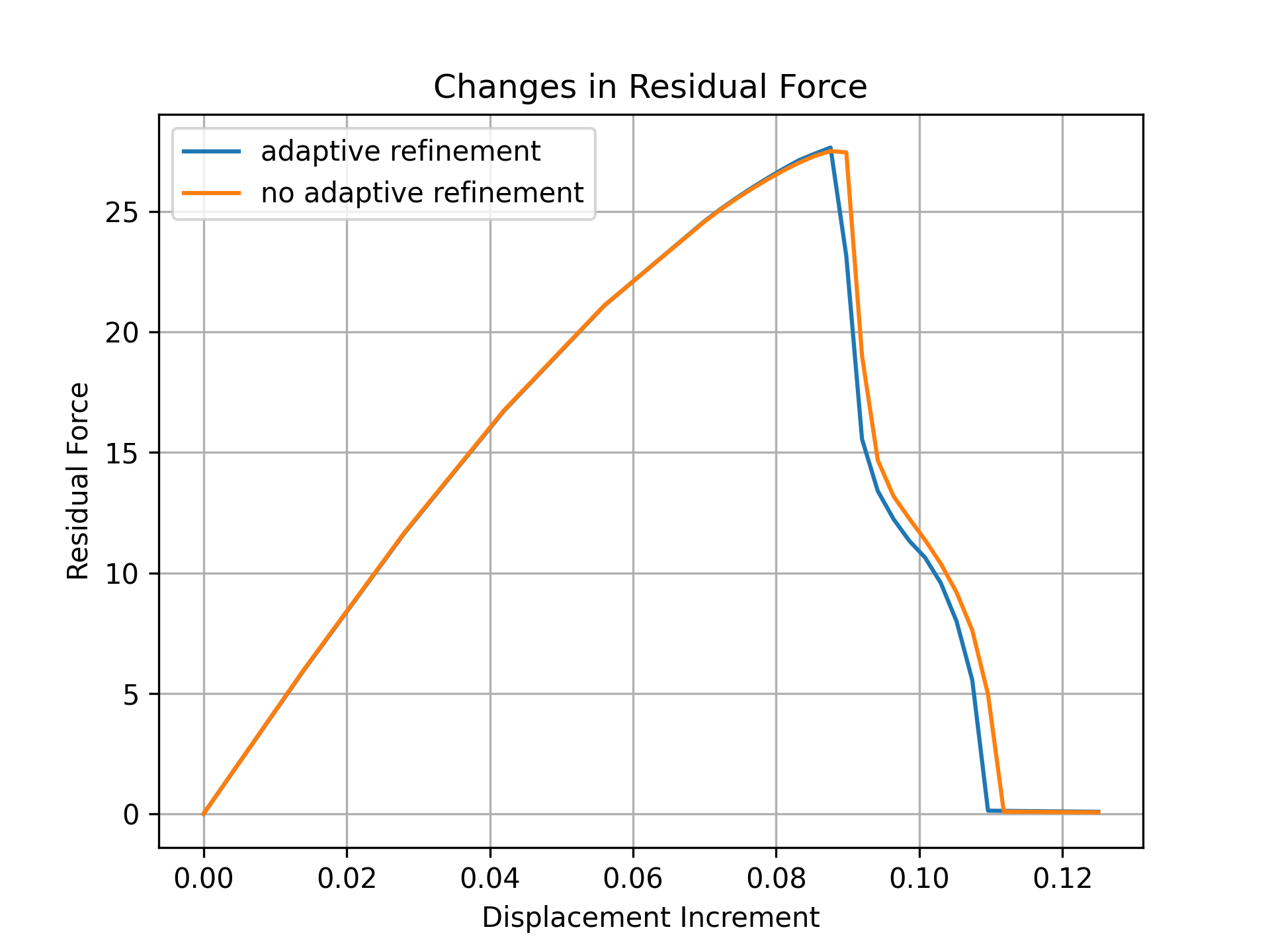}
    \caption{Residual froce with load-displacement curves of the square model with a circular notch,
    adaptive and non-adaptive}
    \label{fig:energy_curves}
\end{figure}

Figure \ref{fig:energy_curves} demonstrates the evolution of residual force
throughout the fracture process under an external load. During the initial
loading of displacement, the residual force gradually increases. Upon the
occurrence of fracture, the stored potential energy begins to convert into
fracture energy, resulting in a gradual decrease in the residual force until the
material completely breaks, at which point the residual force becomes zero. The
trends in the simulations, both with and without adaptive mesh refinement,
remain consistent.  This result is consistent with the result obtained by
Bourdin\cite{Bourdin00}.

The adaptive mesh refinement approach enables a more targeted refinement of the
mesh in regions of interest, such as around crack tips, resulting in improved
accuracy in capturing the energy evolution. In summary, Figure
\ref{fig:energy_curves} highlights the role of adaptive mesh refinement in
improving the accuracy of simulations while maintaining consistency in the
observed residual force changes during fracture.

\subsection{Model with a Notch: Upward Stretch on the Upper Boundary in the $y$ Direction}
\begin{figure}[htbp]
    \centering
    \includegraphics[width=0.3\textwidth]{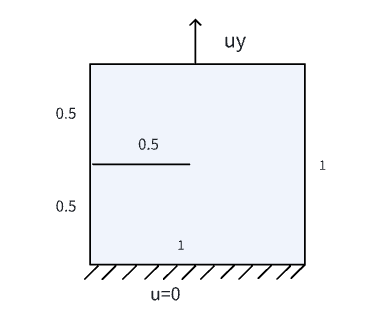}
    \includegraphics[width=0.6\textwidth]{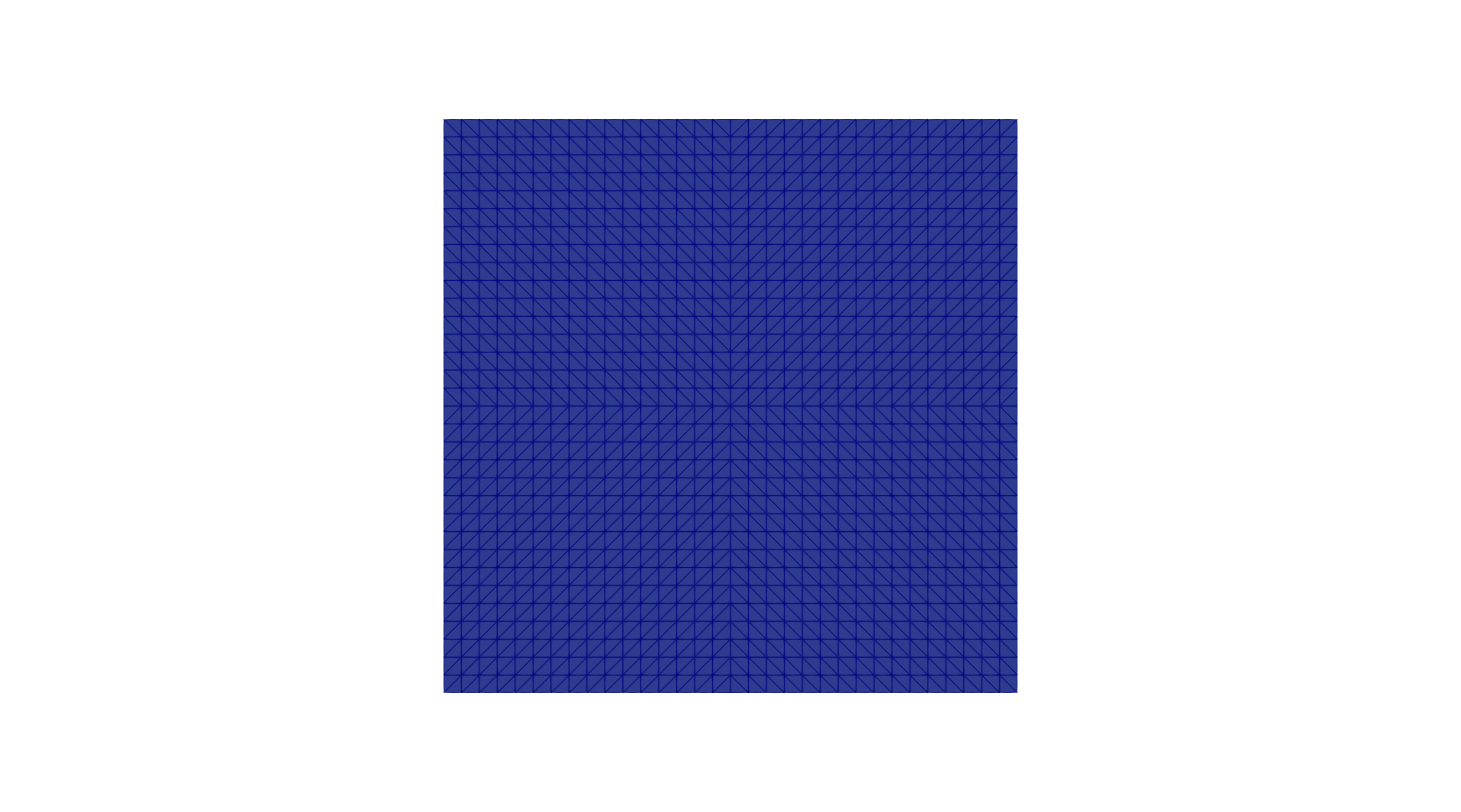}
    \caption{Geometry and boundary conditions of the model with a notch (the
    unit is mm) and initial mesh generation of the model with a notch}
    \label{fig:model_notch1}
\end{figure}

Consider a square model with a unilateral incision, as shown in Figure
\ref{fig:model_notch1}. The model has the following material parameters: the
critical energy release rate is given by $G_c = 2.7 \times 10^{-3} \text{
kN/mm}$, the scale factor $l_0 = 1.33 \times 10^{-2} \text{ mm}$, Lamé's first
parameter $\lambda = 121.15 \text{ kN/mm}^{-2}$, and Lamé's second parameter
$\mu = 80.77 \text{ kN/mm}^{-2}$\cite{Tian19}.

The material's bottom boundary at $y = 0$ is fixed in place. The upper boundary
at $y = 1$ experiences a displacement increment of $\Delta u_y = 1 \times
10^{-5} \text{ mm}$ for the initial 500 steps, followed by an increment of
$\Delta u_y = 1 \times 10^{-6} \text{ mm}$ for subsequent steps, with a total
displacement of $u_y = 6.1 \times 10^{-3} \text{ mm}$.

\begin{figure}[htbp]
    \centering
    \includegraphics[width=\textwidth]{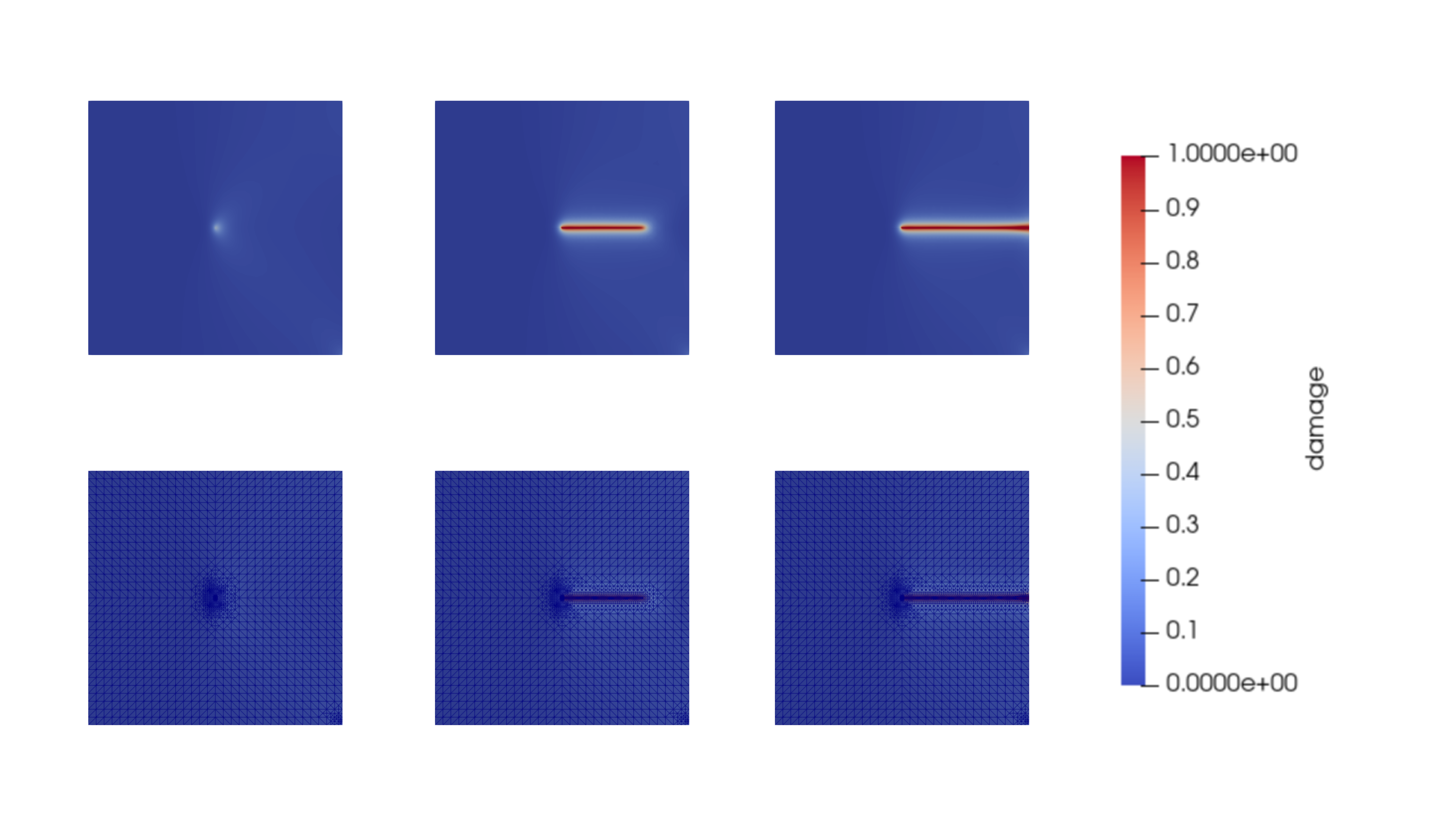}
    \caption{The phase field value (upper) and the mesh refinement
        situation (lower) of the model with a notch (stretch), under different
        displacement loading conditions, from left to right $\Delta u_y$ in
        order $0.0055 \text{ mm}$, $0.0059 \text{ mm}$, $0.006 \text{ mm}$}
    \label{fig:mesh_refinement1}
\end{figure}

Figure \ref{fig:mesh_refinement1} illustrates the adaptive refinement of the
mesh during the crack growth process at three distinct stages. Initially, the
model starts with a coarse mesh consisting of 2048 elements. As the load is
incrementally applied, the mesh undergoes adaptive refinement to accurately
capture the evolving crack geometry. If we don not use adaptive refinement, the
number of elements should not be less than 45000 to ensure the accuracy of the
simulation. So the adaptive refinement method significantly reduces the number
of elements and improves the computational efficiency.

At the first stage, with an incremental displacement of $\Delta u_y =0.0055 mm$
applied, the mesh elements increase to 3406. This refinement is concentrated
around the crack tip and along the expected crack path, indicating that the
simulation is anticipating the direction of crack propagation.

At the second stage, with an additional displacement of $\Delta u_y=0.0059 mm$,
the mesh elements further increase to 6040. The refinement continues to focus on
the crack tip and the surrounding region, ensuring that the simulation captures
the fine details of the crack growth.

Finally, with a further displacement of $\Delta u_y =0.006 mm$, the mesh
elements expand to 7510. The adaptive refinement continues to evolve,
maintaining a high resolution in the critical regions while minimizing the
number of elements in less important areas.

\begin{figure}[htbp]
    \centering
    \includegraphics[width=0.5\textwidth]{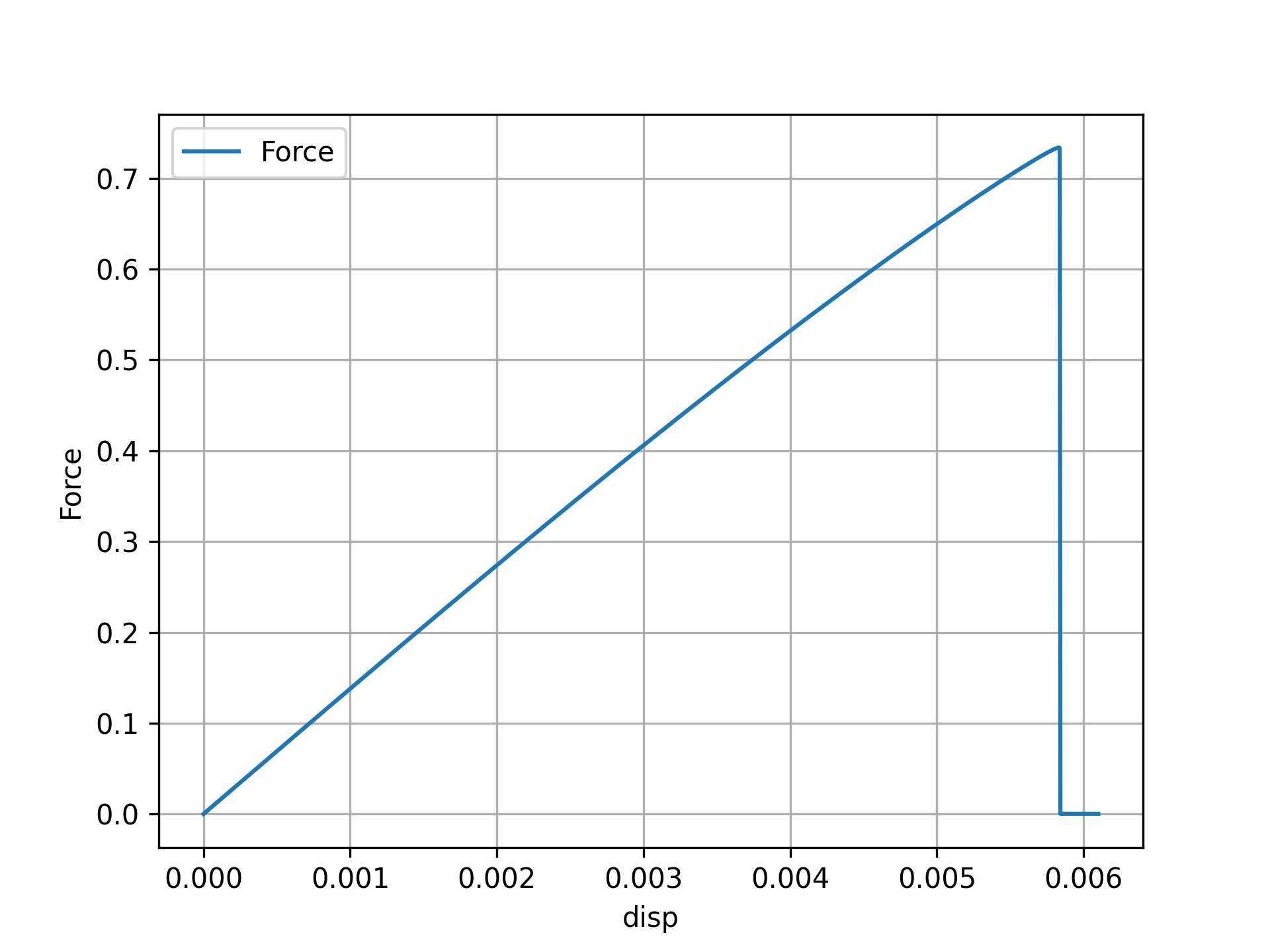}
    \caption{Load-displacement curves in the model with a notch (stretch)}
    \label{fig:load_displacement1}
\end{figure}

Figure \ref{fig:load_displacement1} depicts the evolution of residual force at
the node where the load is applied during the loading process.  Notably, the
observed outcome closely aligns with the previous findings reported by
Miehe\cite{Miehe10a, Miehe10b} and Ambati et al. using conventional phase field
model (PFM) calculations. This consistency strengthens the credibility of the
current results and highlights the reliability of the adaptive finite element
simulation approach.

\subsection{Model with a Notch: Shearing Stretch on the Upper Boundary in the x-Direction}
The model parameters are consistent with those used in the previous model, with
the main change being the modification of the direction of displacement
increment. Specifically, while retaining all other boundary conditions and model
parameters, replace the shear force on the right side with the upward
tension in the above model to verify the accuracy of simulating type II
fracture\cite{Tian19}.

\begin{figure}[htbp]
    \centering
    \includegraphics[width=0.5\textwidth]{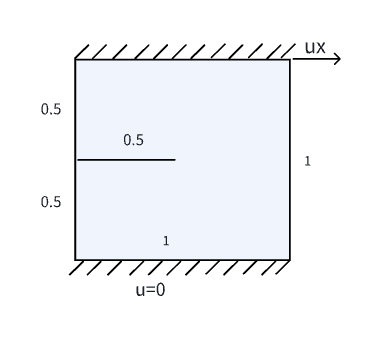}
    \caption{Geometry and boundary conditions of the model with a notch (the unit is mm)}
    \label{fig:notch_x_direction2}
\end{figure}

At each step of the upper boundary where $y = 1$, the displacement increment is
$\Delta u_x = 1 \times 10^{-5} \text{ mm}$, and it loads 1700 steps. The total
displacement increment is $u_x = 1.7 \times 10^{-2} \text{ mm}$. The bottom of
the material is fixed at $y = 0$. Figure \ref{fig:notch_x_direction2} shows
the geometric shape of the model and the initial mesh generation.

\begin{figure}[htbp]
    \centering
    \includegraphics[width=\textwidth]{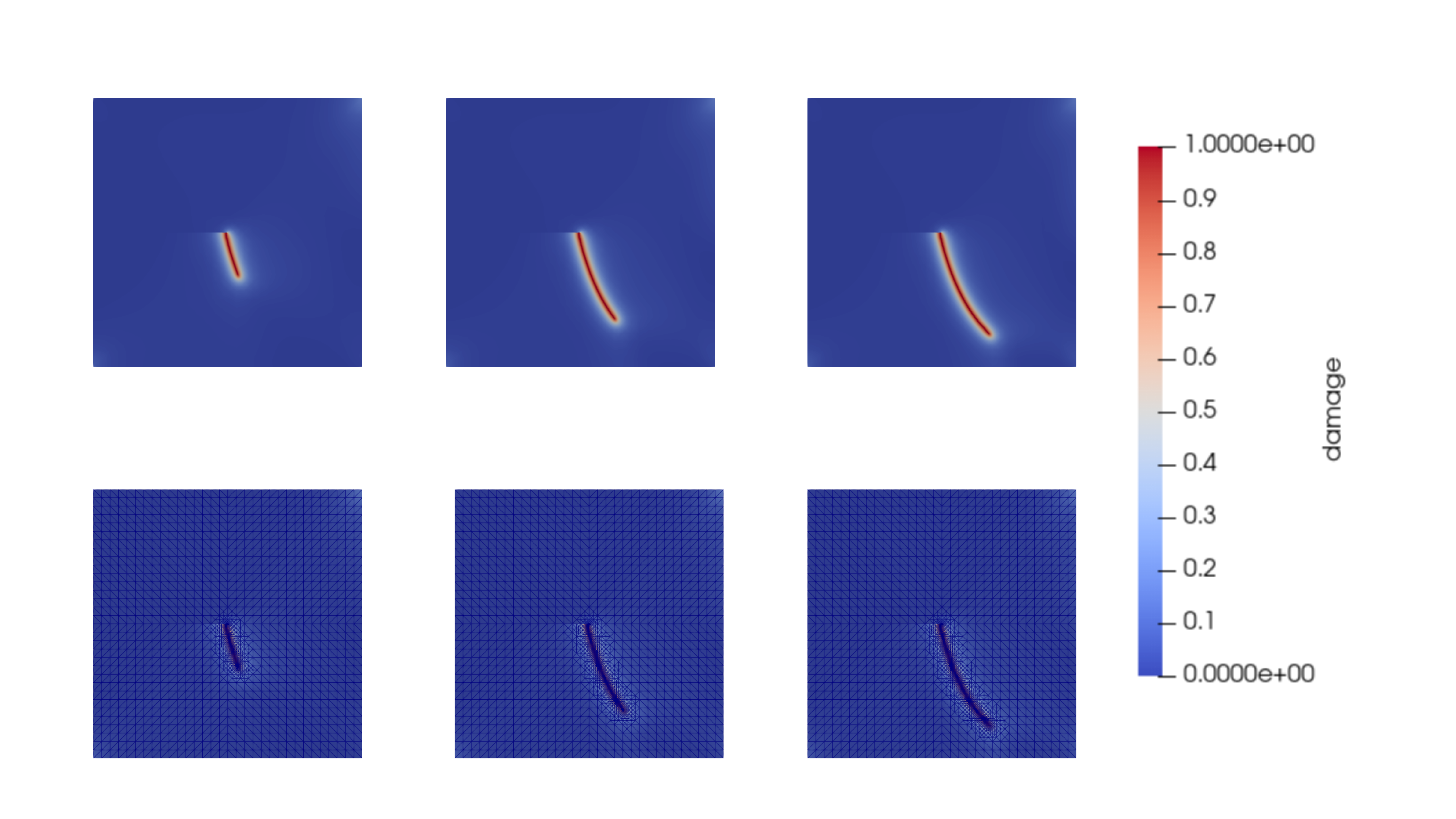}
    \caption{The phase field value (upper) and the mesh refinement
        situation (lower) of the model with a notch (shearing), under different
        displacement loading conditions, from left to right $\Delta u_x$ in
        order $0.012 \text{ mm}, 0.015 \text{ mm}, 0.017 \text{ mm}, 0.021 \text{ mm}$}
    \label{fig:mesh_refinement2}
\end{figure}

Figure \ref{fig:mesh_refinement2} shows the mesh refinement patterns
and phase field values for varying displacement loadings in the x-direction. As
the load increases ($\Delta u_x =$ $0.012 \text{ mm},$ $ 0.015 \text{ mm},$ $ 0.017
\text{ mm},$ $ 0.021 \text{ mm}$), cracks gradually form, and the mesh elements
transition from 2048 to 3786, 5540, 6226, and 7098. Notably, the mesh refinement
consistently aligns with and outpaces crack generation, enabling early capture
of the crack propagation direction.

\begin{figure}[htbp]
    \centering
    \includegraphics[width=0.5\textwidth]{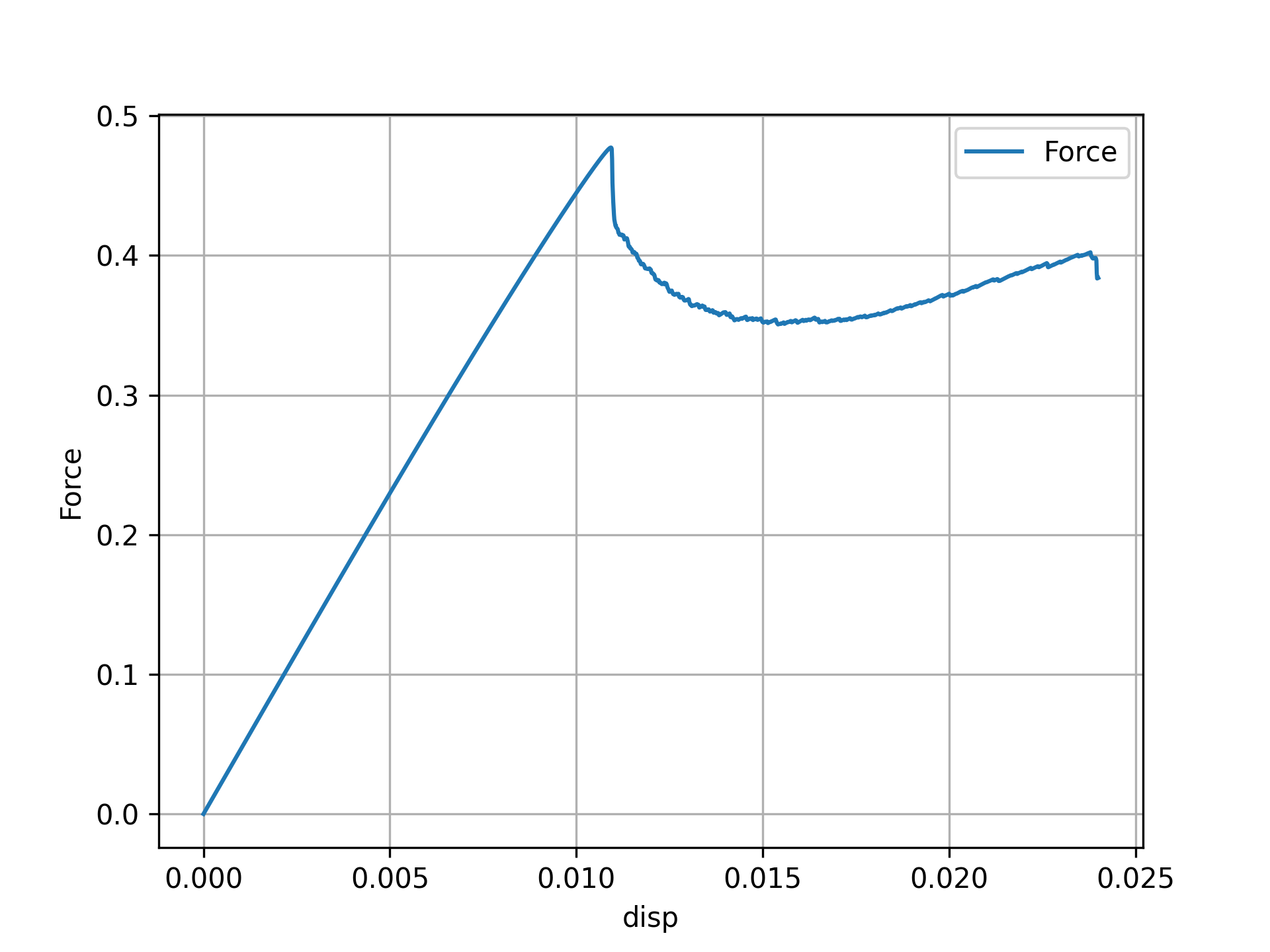}
    \caption{Load-displacement curves in the model with a notch (shearing)}
    \label{fig:load_displacement_curves2}
\end{figure}

Figure \ref{fig:load_displacement_curves2} illustrates the load-displacement
curves for the unilateral incision model. The variation of residual force in the
stressed region upon load application exhibits consistency with results obtained
by Miehe and Ambati using standard phase field model calculations.

These computational results demonstrate that the algorithm effectively simulates
both type I and type II fractures through adaptive mesh refinement, enhancing
computational efficiency without compromising accuracy. The alignment of mesh
refinement with crack propagation further validates the algorithm's ability to
accurately capture fracture behavior.

\subsection{The L-shaped plate (Stretch and compression mixing)}
\begin{figure}[htbp]
    \centering
    \includegraphics[width=0.3\textwidth]{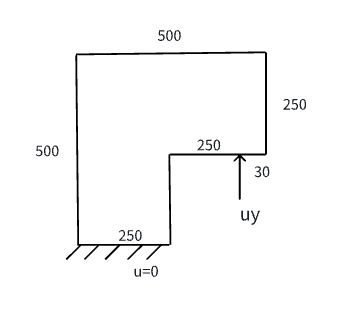}
    \includegraphics[width=0.5\textwidth]{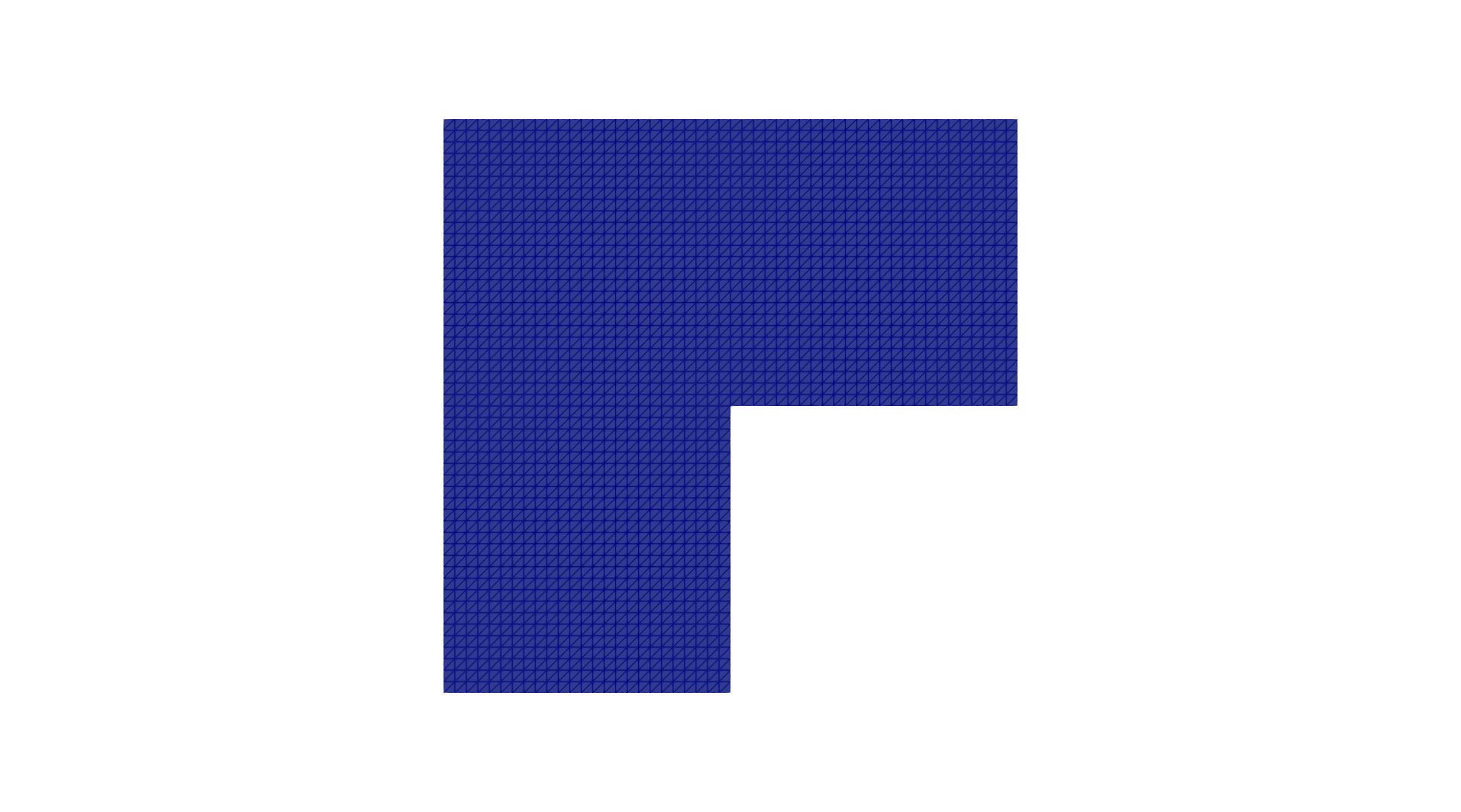}
    \caption{Geometry and boundary conditions of the L-shaped plate (the unit is
    mm) (left), and initial mesh generation of the model (right)}
    \label{fig:l_shaped}
\end{figure}

Consider an L-shaped region, the model has the following material parameters.
The critical energy release rate is given by $G_c=8.9e-5 \text{ kN/mm}$, the
scale factor $l_0=1.88 \text{ mm}$, the Lamé's first parameter
$\lambda=6.16\text{ kN/mm}^{-2}$, and the Lamé's second parameter
$\mu=10.95\text{ kN/mm}^{-2}$. Figure \ref{fig:l_shaped} shows the geometric shape of the model
and the initial mesh generation\cite{Tian19}.

\begin{figure}[htbp]
    \centering
    \includegraphics[width=0.5\textwidth]{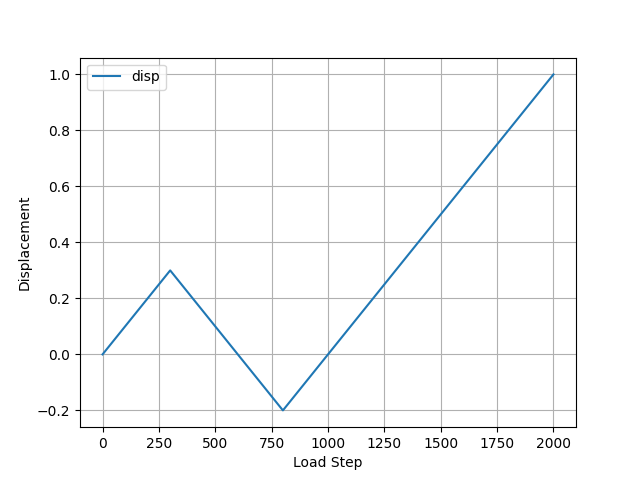}
    \caption{Displacement loading situation of the L-shaped plate}
    \label{fig:l_shaped_plate}
\end{figure}

For this model, we impose a tension-compression mixed load condition at a point,
with specific parameters as illustrated in Figure\ref{fig:l_shaped_plate}.
Initially, an upward load is applied, followed by a downward load, and then an
upward load again. The purpose of this loading sequence is to assess the
algorithm's capability to handle compression crack closure under
tension-compression mixed conditions.

\begin{figure}[htbp]
    \centering
    \includegraphics[width=\textwidth]{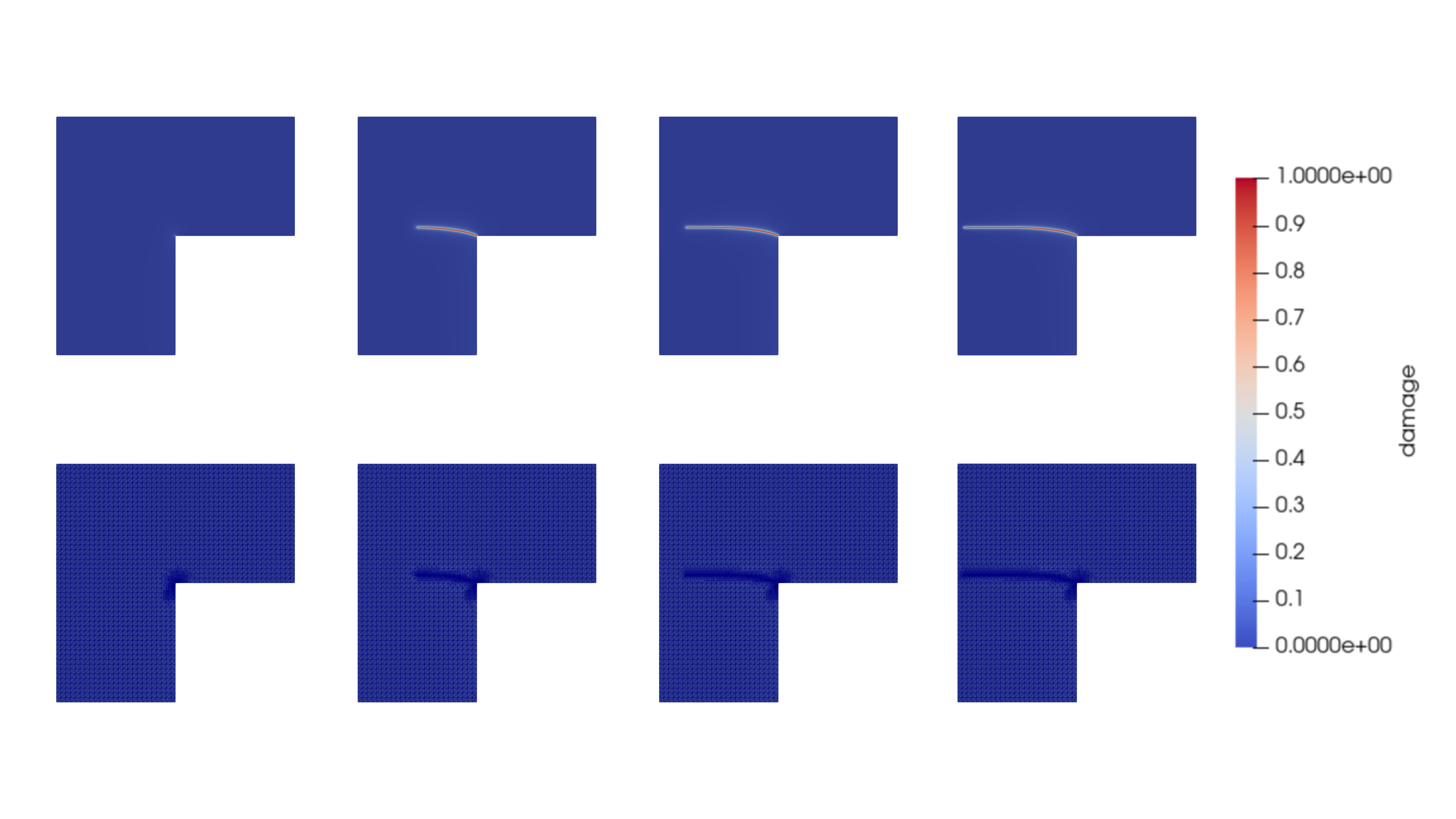}
    \caption{The phase field value (upper) and the mesh refinement
    situation (lower) of the L-shaped plate, under different displacement
    loading conditions, from left to right $\Delta u_y$ in order 0.22\text{ mm},
    0.3\text{ mm}, 0.45\text{ mm}, 1\text{ mm}}
    \label{fig:mesh_refinement_l_shaped}
\end{figure}

In Figure \ref{fig:mesh_refinement_l_shaped}, we observe significant changes in
the mesh refinement patterns and corresponding phase field values as
displacement loading conditions vary. These changes are presented sequentially
from left to right, corresponding to different displacement values of $\Delta
u_y=0.22\text{ mm}, 0.3\text{ mm}, 0.45\text{ mm}, 1\text{ mm}$. As the load is
gradually applied, cracks form progressively in the simulated L-shaped plate.

It is worth noting that the number of mesh elements also changes during the
process of crack formation and propagation. Starting with an initial count of
3750 mesh elements, the number increases to 6563, 22549, 31265, and ultimately
reaches 34419 as cracks emerge and propagate. Under compression conditions,
once a crack forms, it does not close, reflecting the irreversibility of crack
propagation. This characteristic is crucial for understanding and predicting the
failure behavior of materials under compressive loads. 

\begin{figure}[htbp]
    \centering
    \includegraphics[width=0.5\textwidth]{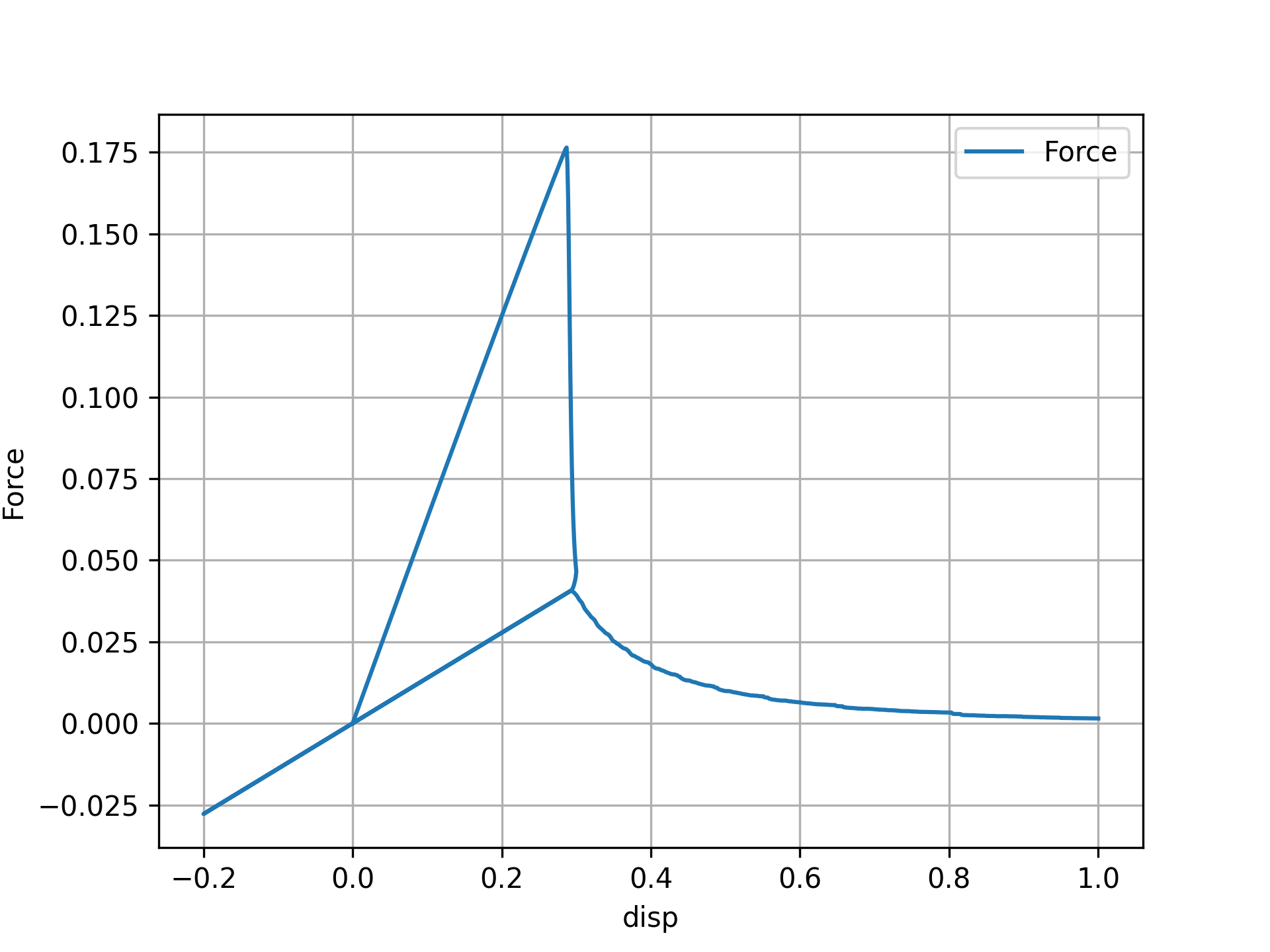}
    \caption{Load-displacement curves of the L-shaped plate}
    \label{fig:load_displacement_curves3}
\end{figure}

Additionally, Figure \ref{fig:load_displacement_curves3} illustrates the
variation of residual force with the application of load, showing excellent
agreement with results calculated by Singh et al.

\subsection{3D model with planar cutouts: Upward Stretch on the Upper Boundary
in the $z$ Direction}

For a 3D model featuring planar slices, the model parameters are set as follows:
a critical energy release rate of $G_c=5e-4 \text{ kN/mm}$, a scale factor of
$l_0=0.2 \text{ mm}$, Young's modulus $E=20.8 \text{ KN/mm}^2$, and Poisson's
ratio $\nu=0.3$\cite{LI201944}.

\begin{figure}[htbp]
    \centering
    \includegraphics[width=0.15\textwidth]{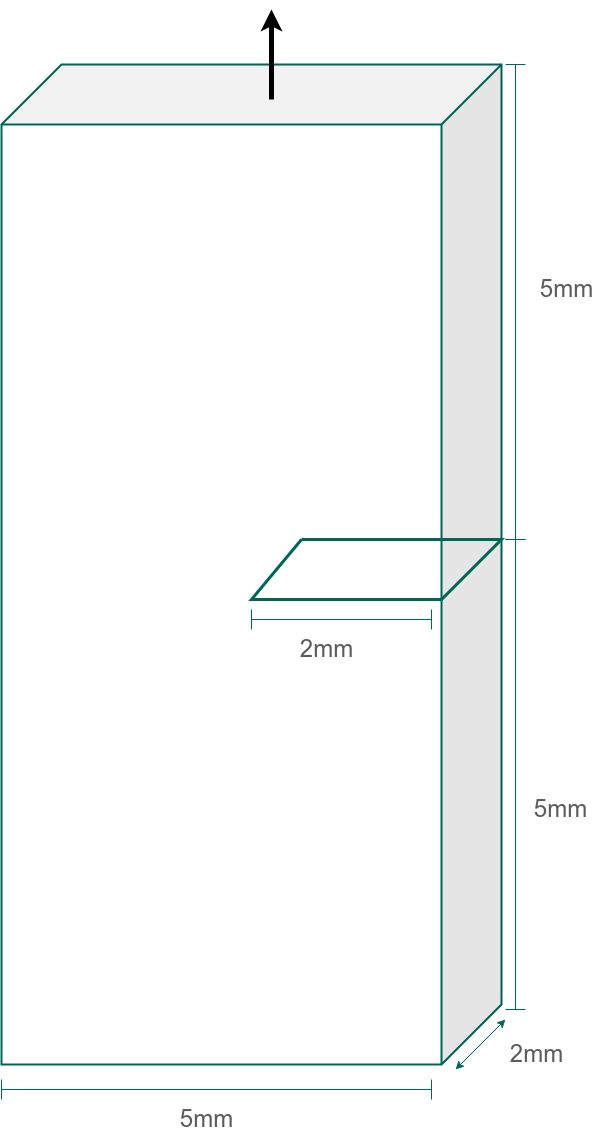}
    \includegraphics[width=0.5\textwidth]{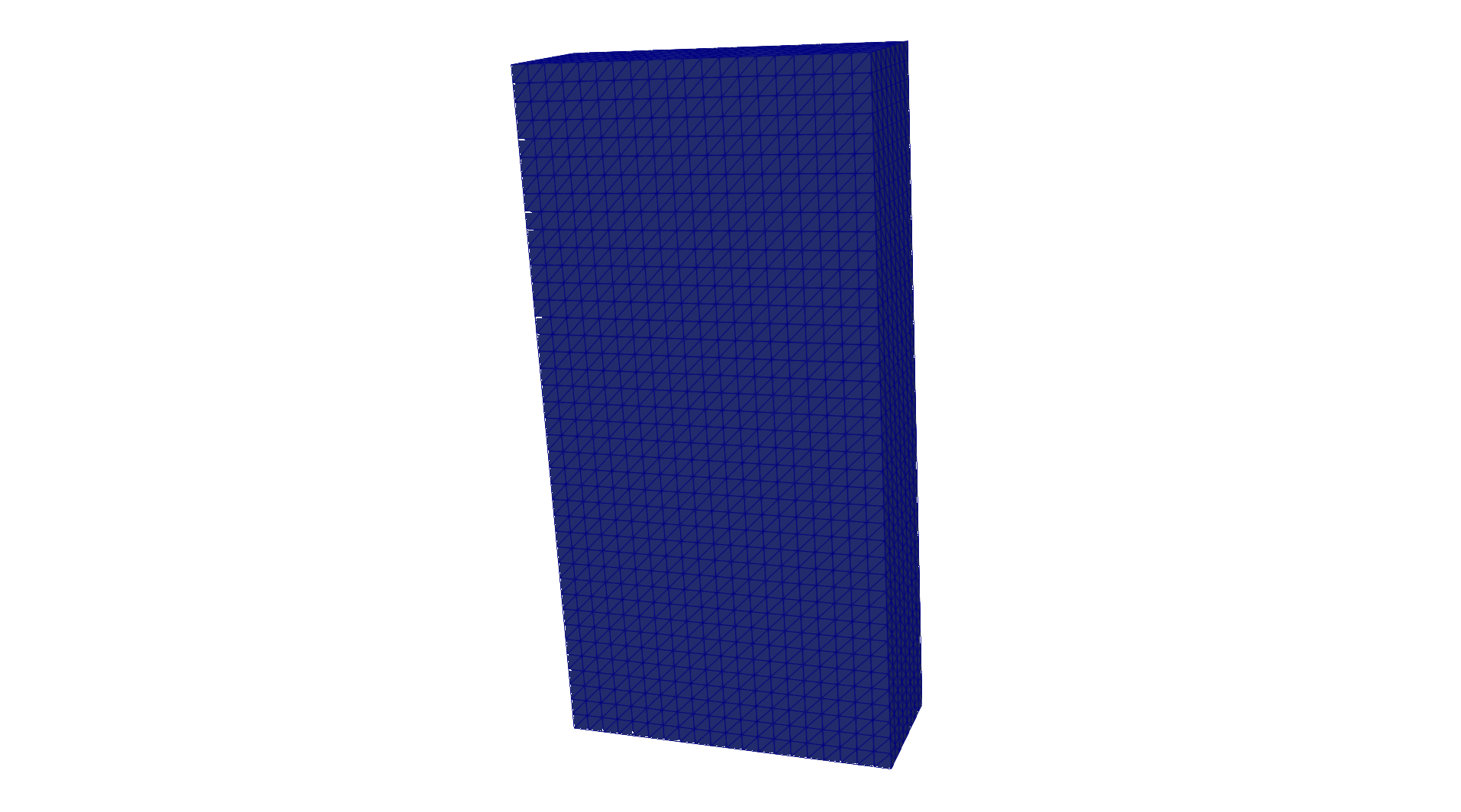}
    \caption{Geometry and boundary conditions (the unit is
    mm) (left), and initial mesh generation of the 3D model (right)}
    \label{fig:initial_mesh_3d}
\end{figure}

The boundary conditions in this model resemble those of a 2D extruded model. An
upward displacement is applied at the upper boundary located at $z=10$, with an
increment of $\Delta u_z=1 \times 10^{-4} \text{ mm}$ resulting in a total
displacement increment of $u_z=4.5 \times 10^{-2} \text{ mm}$. The bottom of the
material is fixed at $z=0$. Figure \ref{fig:initial_mesh_3d} shows the geometric
shape of the model and the initial mesh generation.

In the three-dimensional example, we employed a tetrahedral mesh for initial
mesh generation, utilizing the bisection method for mesh refinement during the
adaptive process.

\begin{figure}[htbp]
    \centering
    \includegraphics[width=\textwidth]{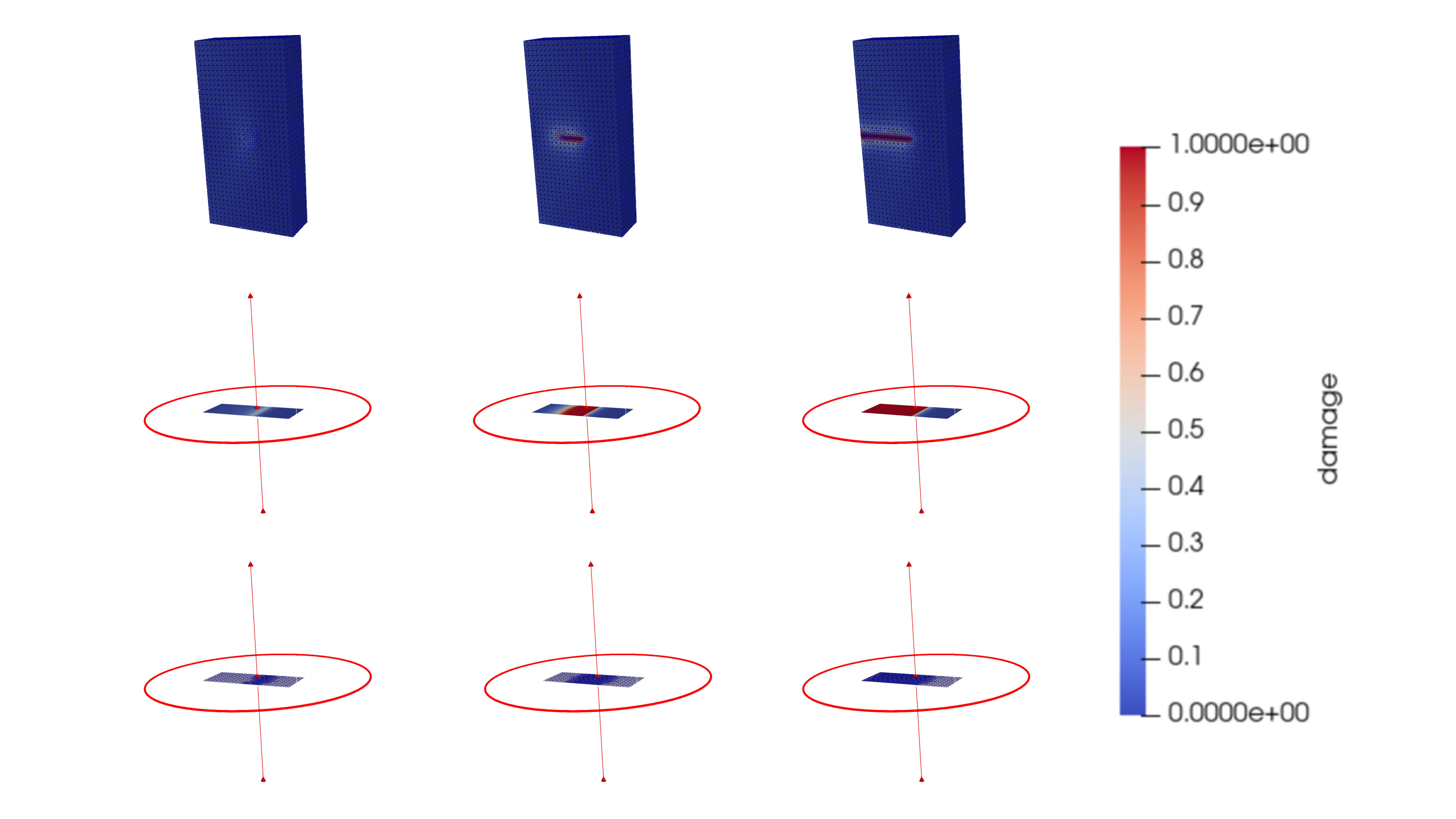}
    \caption{The phase field value (upper) and the mesh refinement
    situation (lower) of the 3D model, under different displacement loading
    conditions, from left to right $\Delta u_z$ in order 0.27\text{ mm}, 0.28\text{
    mm}, 0.29\text{ mm}}
    \label{fig:mesh_refinement_3d}
\end{figure}

\begin{figure}[htbp]
    \centering
    \includegraphics[width=0.5\textwidth]{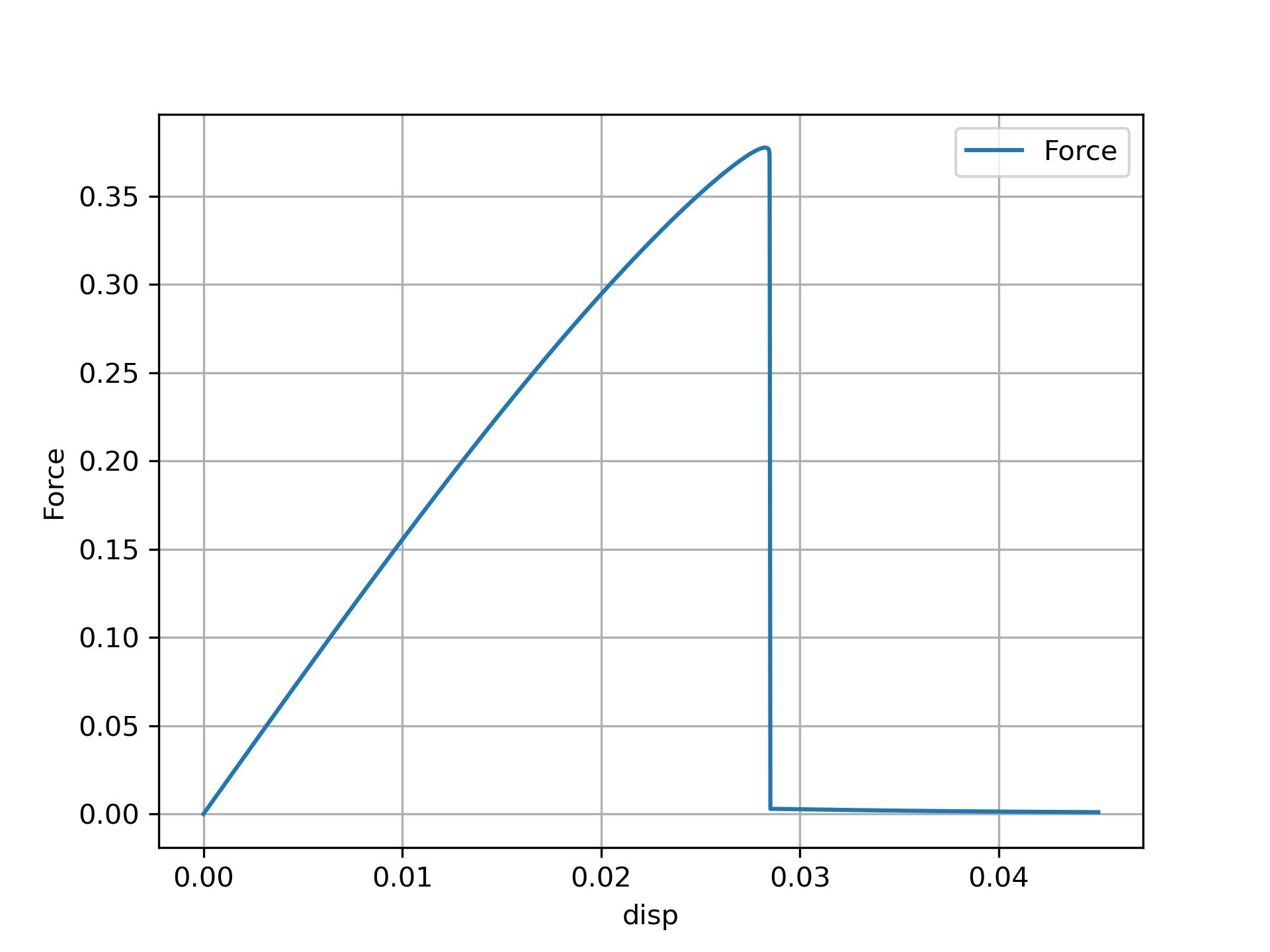}
    \caption{Load-displacement curves of the 3D model}
    \label{fig:load_displacement_curves_3d}
\end{figure}

By examining the variation diagram of residual force, we can observe that the
results obtained from the three-dimensional simulation are in good agreement
with those of the corresponding two-dimensional calculations. This consistency
validates the applicability and reliability of the chosen numerical methods in
both two- and three-dimensional settings. Furthermore, the trend of mesh
refinement observed in the three-dimensional simulation aligns closely with the
direction of crack development. As the load is gradually applied, cracks begin
to form, and the number of mesh elements increases accordingly from the initial
38400 to 52352, 64328, and ultimately 82289. This progressive refinement ensures
that the mesh maintains adequate resolution in areas where crack propagation is
occurring, thereby improving the accuracy of the simulation results.

\section{Summary}
\label{sec:summary}

This paper presents an adaptive finite element method (AFEM) for phase field
fracture models, leveraging recovery-type posterior error estimates to achieve
efficient mesh refinement. The method transforms the gradient of the numerical
phase field into a smoother function space, using the discrepancy between the
recovered and original gradients as an error indicator. This approach eliminates
the need for empirical parameters, allowing the automatic and accurate capture
of crack propagation directions while significantly reducing the number of mesh
elements and improving computational efficiency. By integrating these
innovations into the FEALPy platform, we have implemented a robust phase field
fracture simulation module with flexible refinement techniques and
GPU-accelerated matrix solvers. Numerical experiments on classical 2D and 3D
brittle fracture problems validate the accuracy, robustness, and efficiency of
the proposed method.  

In future work, we plan to expand this framework by
developing a dedicated fracture numerical simulation application (APP) built on
FEALPy. This APP will serve as a comprehensive platform for the integration,
validation, and comparison of various fracture simulation algorithms, providing
a rich library of computational examples. Leveraging FEALPy’s modular
architecture, multi-backend tensor computation engine, and extensive algorithmic
support, the APP will offer diverse simulation capabilities, supporting scalable
tensor computation with backend options such as Numpy, PyTorch, and JAX. The
application will harness FEALPy's powerful functionality, including hierarchical
and modular structures, and an array of core algorithms, to facilitate seamless,
high-performance simulations on heterogeneous hardware systems. Additionally,
the APP will enable the generation of large-scale datasets for machine learning
models, fostering interdisciplinary research in fracture mechanics and
intelligent simulation technologies.


\section*{Acknowledgments}
This work were supported by the National Natural Science Foundation of China
(NSFC) (Grant Nos. 12371410, 12261131501), the construction of innovative
provinces in Hunan Province (Grant No. 2021GK1010), and the Graduate Innovation
Project of Xiangtan University (Nos. XDCX2023Y135, XDCX2024Y175).

During the preparation of this work the author(s) used ChatGPT 4o in order to
improve language and readability. After using this tool/service, the author(s)
reviewed and edited the content as needed and take(s) full responsibility for
the content of the publication.

\bibliographystyle{unsrtnat}

\end{document}